# QUANTUM COHOMOLOGY
# OF PROJECTIVE BUNDLES OVER $\mathbb{P}^n$

Zhenbo Qin[1] and Yongbin Ruan[2]

## 1. Introduction

Quantum cohomology, proposed by Witten's study [16] of two dimensional nonlinear sigma models, plays a fundamental role in understanding the phenomenon of mirror symmetry for Calabi-Yau manifolds. This phenomenon was first observed by physicists motivated by topological field theory. A topological field theory starts with correlation functions. The correlation functions of sigma model are linked with the intersection numbers of cycles in the moduli space of holomorphic maps from Riemann surfaces to manifolds. For some years, the mathematical construction of these correlation functions remained to be a difficult problem because the moduli spaces of holomorphic maps usually are not compact and may have wrong dimension. The quantum cohomology theory was first put on a firm mathematical footing by [12,13] for semi-positive symplectic manifolds (including Fano and Calabi-Yau manifolds), using the method of symplectic topology. Recently, an algebro-geometric approach has been taken by [8,9]. The results of [12,13] have been redone in the algebraic geometric setting for the case of homogeneous spaces. The advantage of homogeneous spaces is that the moduli spaces of holomorphic maps always have expected dimension and their compactifications are nice. Beyond the homogeneous spaces, one can not expect such nice properties for the moduli spaces. The projective bundles are perhaps the simplest examples. However, by developing sophisticated excessive intersection theory, it is possible that the algebro-geometric method can work for any projective manifolds. In turn, it may shed new light to removing the semi-positive condition in the symplectic setting.

Although we have a solid foundation for quantum cohomology theory at least for semi-positive symplectic manifolds, the calculation remains to be a difficult task. So far, there are only a few examples which have been computed, e.g., Grassmannian [14], some rational surfaces [6], flag varieties [4], some complete intersections [3], and the moduli space of stable bundles over Riemann surfaces [15]. One of the common feature for these examples is that the relevant moduli spaces of rational curves have expected dimension. Then, one can use the intersection theory. We should mention that there are many predications based on mathematically unjustified mirror symmetry (for Calabi-Yau 3-folds) and linear sigma model (for toric varieties). In this paper, we attempt to determine the quantum cohomology of

[1] Partially supported by an NSF grant
[2] Partially supported by an NSF grant and a Sloan fellowship

Typeset by $\mathcal{A}_{\mathcal{M}}\mathcal{S}$-TeX

projective bundles over the projective space $\mathbb{P}^n$. In contrast to the previous examples, the relevant moduli spaces in our case frequently do not have expected dimensions. It makes the calculation more difficult. We overcome this difficulty by using excessive intersection theory.

There are two main ingredients in our arguments. The first one is a result of Siebert and Tian (the Theorem 2.2 in [14]), which says that if the ordinary cohomology $H^*(X;\mathbb{Z})$ of a symplectic manifold $X$ with the symplectic form $\omega$ is the ring generated by $\alpha_1,\ldots,\alpha_s$ with the relations $f^1,\ldots,f^t$, then the quantum cohomology $H^*_\omega(X;\mathbb{Z})$ of $X$ is the ring generated by $\alpha_1,\ldots,\alpha_s$ with $t$ new relations $f^1_\omega,\ldots,f^t_\omega$ where each new relation $f^i_\omega$ is just the relation $f^i$ evaluated in the quantum cohomology ring structure. It was known that the quantum product $\alpha \cdot \beta$ is the deformation of ordinary cup product by the lower order terms called quantum corrections. The second ingredient is that under certain numerical conditions, most of the quantum corrections vanishes. Moreover, the nontrivial quantum corrections seem to come from Mori's extremal rays.

Let $V$ be a rank-$r$ bundle over $\mathbb{P}^n$, and $\mathbb{P}(V)$ be the corresponding projective bundle. Let $h$ and $\xi$ be the cohomology classes of a hyperplane in $\mathbb{P}^n$ and the tautological line bundle in $\mathbb{P}(V)$ respectively. For simplicity, we make no distinction between $h$ and $\pi^*h$ where $\pi : \mathbb{P}(V) \to \mathbb{P}^n$ is the natural projection. Denote the product of $i$ copies of $h$ and $j$ copies of $\xi$ in the ordinary cohomology ring by $h_i\xi_j$, and the product of $i$ copies of $h$ and $j$ copies of $\xi$ in the quantum cohomology ring by $h^i \cdot \xi^j$. For $i = 0,\ldots,r$, put $c_i(V) = c_i \cdot h_i$ for some integer $c_i$. It is well known that $-K_{\mathbb{P}(V)} = (n+1-c_1)h + r\xi$ and the ordinary cohomology ring $H^*(\mathbb{P}(V);\mathbb{Z})$ is the ring generated by $h$ and $\xi$ with the two relations:

$$h_{n+1} = 0 \qquad \text{and} \qquad \sum_{i=0}^{r}(-1)^i c_i \cdot h_i \xi_{r-i} = 0. \tag{1.1}$$

In particular, $H^{2(n+r-2)}(\mathbb{P}(V);\mathbb{Z})$ is generated by $h_{n-1}\xi_{r-1}$ and $h_n\xi_{r-2}$, and its Poincaré dual $H_2(\mathbb{P}(V);\mathbb{Z})$ is generated by $(h_{n-1}\xi_{r-1})_*$ and $(h_n\xi_{r-2})_*$ where for $\alpha \in H^*(\mathbb{P}(V);\mathbb{Z})$, $\alpha_*$ stands for its Poincaré dual. We have

$$-K_{\mathbb{P}(V)}(A) = a(n+1-c_1) + r \cdot \xi(A) = a(n+1-c_1) + r(ac_1 + b) \tag{1.2}$$

for $A = (ah_{n-1}\xi_{r-1} + bh_n\xi_{r-2})_* \in H_2(\mathbb{P}(V);\mathbb{Z})$.

By definition, $V$ is an ample (respectively, nef) bundle if and only if the tautological class $\xi$ is an ample (respectively, nef) divisor on $\mathbb{P}(V)$. Assume that $V$ is ample such that either $c_1 \leq (n+1)$ or $c_1 \leq (n+r)$ and $V \otimes \mathcal{O}_{\mathbb{P}^n}(-1)$ is nef. Then both $\xi$ and $-K_{\mathbb{P}(V)}$ are ample divisors. Thus, $\mathbb{P}(V)$ is a Fano variety, and its quantum cohomology ring is well-defined [13]. Here we choose the symplectic form $\omega$ on $\mathbb{P}(V)$ to be the Kahler form $\omega$ such that $[\omega] = -K_{\mathbb{P}(V)}$. Let $f^1_\omega$ and $f^2_\omega$ be the two relations in (1.1) evaluated in the quantum cohomology ring $H^*_\omega(\mathbb{P}(V);\mathbb{Z})$. Then by the Theorem 2.2 in [14], the quantum cohomology $H^*_\omega(\mathbb{P}(V);\mathbb{Z})$ is the ring generated by $h$ and $\xi$ with the two relations $f^1_\omega$ and $f^2_\omega$:

$$H^*_\omega(\mathbb{P}(V);\mathbb{Z}) = \mathbb{Z}[h,\xi]/(f^1_\omega, f^2_\omega) \tag{1.3}$$

By Mori's Cone Theorem [5], $\mathbb{P}(V)$ has exactly two extremal rays $R_1$ and $R_2$. Up to an order of $R_1$ and $R_2$, the integral generator $A_1$ of $R_1$ is represented by lines in the fibers of the projection $\pi$. We shall show that under certain numerical conditions, the nontrivial homology classes $A \in H_2(\mathbb{P}(V); \mathbb{Z})$ which give nontrivial quantum corrections are $A_1$ and $A_2$, where $A_2$ is represented by some smooth rational curves in $\mathbb{P}(V)$ which are isomorphic to lines in $\mathbb{P}^n$ via $\pi$. In general, it is unclear whether $A_2$ generates the second extremal ray $R_2$. However, we shall prove that under further restrictions on $V$, $A_2$ generates the extremal ray $R_2$. These analyses enable us to determine the quantum cohomology ring $H^*_\omega(\mathbb{P}(V); \mathbb{Z})$.

The simplest ample bundle over $\mathbb{P}^n$ is perhaps the direct sum of line bundles $V = \oplus_{i=1}^r \mathcal{O}(m_i)$ where $m_i > 0$ for every $i$. Since we can twist $V$ by $\mathcal{O}(-1)$ without changing $\mathbb{P}(V)$, we can assume that $\min\{m_1, \ldots, m_r\} = 1$. In this case, $\mathbb{P}(V)$ is a special case of toric variety. Batyrev [2] conjectured a general formula for quantum cohomology of toric varieties. Furthermore, he computed the contributions from certain moduli spaces of holomorphic maps which have expected dimensions. In our case, the contributions Batyrev computed are only part of the data to compute the quantum cohomology. As we explained earlier, the difficulty in our case lies precisely in computing the contributions from the moduli spaces with wrong dimensions. Nevertheless, in our case, Batyrev's formula (see also [1]) reads as follows.

**Batyrev's Conjecture:** *Let $V = \oplus_{i=1}^r \mathcal{O}(m_i)$ where $m_i > 0$ for every $i$. Then the quantum cohomology ring $H^*_\omega(\mathbb{P}(V); \mathbb{Z})$ is generated by $h$ and $\xi$ with two relations*

$$h^{n+1} = \prod_{i=1}^r (\xi - m_i h)^{m_i - 1} \cdot e^{-t(n+1+r-\sum_{i=1}^r m_i)} \qquad and \qquad \prod_{i=1}^r (\xi - m_i h) = e^{-tr}.$$

Our first result partially verifies Batyrev's conjecture.

**Theorem A.** *Batyrev's conjecture holds if*

$$\sum_{i=1}^r m_i < \min(2r, (n+1+2r)/2, (2n+2+r)/2).$$

Note that under the numerical condition of Theorem A, only extremal rational curves with fundamental classes $A_1$ and $A_2$ give the contributions to the two relations in the quantum cohomology. The moduli space of rational curves $\mathfrak{M}(A_2, 0)$ with fundamental class $A_2$ does not have expected dimension in general. But it is compact. This fact simplifies a great deal of the excessive intersection theory involved. To remove the numerical condition, we have to consider other moduli spaces (for example $\mathfrak{M}(kA_2, 0)$ with $k > 1$ and its excessive intersection theory). These moduli spaces are not compact in general. Then, we have an extra difficulty of the compactification and the appropriate excessive intersection theory with it. It seems to be a difficult problem and we shall not pursue here.

In general, ample bundles over $\mathbb{P}^n$ are not direct sums of line bundles. We can say much less about its quantum cohomology. However, we obtain some result about its general form and compute the leading coefficient.

**Theorem B.** (i) Let $V$ be a rank-$r$ ample bundle over $\mathbb{P}^n$. Assume either $c_1 \leq n$ or $c_1 \leq (n+r)$ and $V \otimes \mathcal{O}_{\mathbb{P}^n}(-1)$ is nef so that $\mathbb{P}(V)$ is Fano. Then the quantum cohomology $H^*_\omega(\mathbb{P}(V); \mathbb{Z})$ is the ring generated by $h$ and $\xi$ with two relations

$$h^{n+1} = \sum_{i+j \leq (c_1-r)} a_{i,j} \cdot h^i \cdot \xi^j \cdot e^{-t(n+1-i-j)}$$

$$\sum_{i=0}^{r} (-1)^i c_i \cdot h^i \cdot \xi^{r-i} = e^{-tr} + \sum_{i+j \leq (c_1-n-1)} b_{i,j} \cdot h^i \cdot \xi^j \cdot e^{-t(r-i-j)}$$

where the coefficients $a_{i,j}$ and $b_{i,j}$ are integers depending on $V$;

(ii) If we further assume that $c_1 < 2r$, then the leading coefficient $a_{0,c_1-r} = 1$.

It is understood that when $c_1 \leq n$, then the summation $\sum_{i+j \leq (c_1-n-1)}$ in the second relation in Theorem B (i) does not exist. In general, it is not easy to determine all the integers $a_{i,j}$ and $b_{i,j}$ in Theorem B (i). However, it is possible to compute these numbers when $(c_1 - r)$ is relatively small. For instance, when $(c_1 - r) = 0$, then necessarily $V = \mathcal{O}_{\mathbb{P}^n}(1)^{\oplus r}$ and it is well-known that the quantum cohomology $H^*_\omega(\mathbb{P}(V); \mathbb{Z})$ is the ring generated by $h$ and $\xi$ with the two relations $h^{n+1} = e^{-t(n+1)}$ and $\sum_{i=0}^r (-1)^i c_i \cdot h^i \cdot \xi^{r-i} = e^{-tr}$. When $(c_1 - r) = 1$ and $r < n$, then necessarily $V = \mathcal{O}_{\mathbb{P}^n}(1)^{\oplus(r-1)} \oplus \mathcal{O}_{\mathbb{P}^n}(2)$. When $(c_1 - r) = 1$ and $r = n$, then $V = \mathcal{O}_{\mathbb{P}^n}(1)^{\oplus(r-1)} \oplus \mathcal{O}_{\mathbb{P}^n}(2)$ or $V = T_{\mathbb{P}^n}$ the tangent bundle of $\mathbb{P}^n$. In these cases, $V \otimes \mathcal{O}_{\mathbb{P}^n}(-1)$ is nef. In particular, the direct sum cases have been computed by Theorem A. We shall prove the following.

**Proposition C.** The quantum cohomology ring $H^*_\omega(\mathbb{P}(T_{\mathbb{P}^n}); \mathbb{Z})$ with $n \geq 2$ is the ring generated by $h$ and $\xi$ with the two relations:

$$h^{n+1} = \xi \cdot e^{-tn} \qquad \text{and} \qquad \sum_{i=0}^{n} (-1)^i c_i \cdot h^i \cdot \xi^{n-i} = (1 + (-1)^n) \cdot e^{-tn}.$$

Recall that for an arbitrary projective bundle over a general manifold, its cohomology ring is a module over the cohomology ring of the base with the generator $\xi$ and the second relation of (1.1). Naively, one may think that the quantum cohomology of projective bundle is a module over the quantum cohomology of base with the generator $\xi$ and the quantanized second relation. Our calculation shows that one can not expect such simplicity for its quantum cohomology ring. We hope that our results could shed some light on the quantum cohomology for general projective bundles, which we shall leave for future research.

Our paper is organized as follows. In section 2, we discuss the extremal rays and extremal rational curves. In section 3, we review the definition of quantum product and compute some Gromov-Witten invariants. In the remaining three sections, we prove Theorem B, Theorem A, and Proposition C respectively.

**Acknowledgements:** We would like to thank Sheldon Katz, Yungang Ye, and Qi Zhang for valuable helps and stimulating discussions. In particular, we are grateful to Sheldon Katz for bringing us the attention of Batyrev's conjecture.

## 2. Extremal rational curves

Assume that $V$ is ample such that either $c_1 \leq (n+1)$ or $c_1 \leq (n+r)$ and $V \otimes \mathcal{O}_{\mathbb{P}^n}(-1)$ is nef. In this section, we study the extremal rays and extremal rational curves in the Fano variety $\mathbb{P}(V)$. By Mori's Cone Theorem (p.25 in [5]), $\mathbb{P}(V)$ has precisely two extremal rays $R_1 = \mathbb{R}_{\geq 0} \cdot A_1$ and $R_2 = \mathbb{R}_{\geq 0} \cdot A_2$ such that the cone $\mathrm{NE}(\mathbb{P}(V))$ of curves in $\mathbb{P}(V)$ is equal to $R_1 + R_2$ and that $A_1$ and $A_2$ are the homology classes of two rational curves $E_1$ and $E_2$ in $\mathbb{P}(V)$ with $0 < -K_{\mathbb{P}(V)}(A_i) \leq \dim(\mathbb{P}(V))+1$. Up to orders of $A_1$ and $A_2$, we have $A_1 = (h_n \xi_{r-2})_*$, that is, $A_1$ is represented by lines in the fibers of $\pi$. It is also well-known that if $V = \oplus_{i=1}^r \mathcal{O}_{\mathbb{P}^n}(m_i)$ with $m_1 \leq \ldots \leq m_r$, then $A_2 = [h_{n-1}\xi_{r-1} + (m_1-c_1)h_n\xi_{r-2}]_*$ which is represented by a smooth rational curve in $\mathbb{P}(V)$ isomorphic to a line in $\mathbb{P}^n$ via $\pi$. However, in general, it is not easy to determine the homology class $A_2$ and the extremal rational curves representing $A_2$. Assume that

$$V|_\ell = \oplus_{i=1}^r \mathcal{O}_\ell(m_i) \tag{2.1}$$

for generic lines $\ell \subset \mathbb{P}^n$ where we let $m_1 \leq \ldots \leq m_r$. Since $V$ is ample, $m_1 \geq 1$.

**Lemma 2.2.** *Let $A = [h_{n-1}\xi_{r-1} + (m_1 - c_1)h_n\xi_{r-2}]_*$. Then,*
  (i) *$A$ is represented by a smooth rational curve isomorphic to a line in $\mathbb{P}^n$;*
  (ii) *$A_2 = A$ if and only if $(\xi - m_1 h)$ is nef;*
  (iii) *$A_2 = A$ if $2c_1 \leq (n+1)$;*
  (iv) *$A$ can not be represented by reducible or nonreduced curves if $m_1 = 1$.*

*Proof.* (i) Let $\ell \subset \mathbb{P}^n$ be a generic line. Then we have a natural projection $V|_\ell = \oplus_{i=1}^r \mathcal{O}_\ell(m_i) \to \mathcal{O}_\ell(m_1)$. By the Proposition 7.12 in Chapter II of [7], this surjective map $V|_\ell \to \mathcal{O}_\ell(m_1) \to 0$ induces a morphism $g : \ell \to \mathbb{P}(V)$. Then $g(\ell)$ is isomorphic to $\ell$ via the projection $\pi$. Since $h([g(\ell)]) = 1$ and $\xi([g(\ell)]) = m_1$, we have

$$[g(\ell)] = [h_{n-1}\xi_{r-1} + (m_1 - c_1)h_n\xi_{r-2}]_* = A.$$

(ii) First of all, if $A_2 = [h_{n-1}\xi_{r-1}+(m_1-c_1)h_n\xi_{r-2}]_*$, then for any curve $E$, $[E] = a(h_n\xi_{r-2})_* + b[h_{n-1}\xi_{r-1} + (m_1-c_1)h_n\xi_{r-2}]_*$ for some nonnegative numbers $a$ and $b$; so $(\xi - m_1 h)([E]) = a \geq 0$; therefore $(\xi - m_1 h)$ is nef. Conversely, if $(\xi - m_1 h)$ is nef, then $0 \leq (\xi - m_1 h)([E]) = ac_1 + b - am_1$ where $[E] = (ah_{n-1}\xi_{r-1} + bh_n\xi_{r-2})_*$ for some curve $E$; thus $[E] = (ac_1 + b - am_1)(h_n\xi_{r-2})_* + a[h_{n-1}\xi_{r-1} + (m_1 - c_1)h_n\xi_{r-2}]_*$; it follows that $A_2 = [h_{n-1}\xi_{r-1} + (m_1 - c_1)h_n\xi_{r-2}]_* = A$.

(iii) Let $A_2 = (ah_{n-1}\xi_{r-1} + bh_n\xi_{r-2})_*$. Since $A_1 = (h_n\xi_{r-2})_*$ and $a = h(A_2) \geq 0$, $a \geq 1$. If $a > 1$, then since $2c_1 \leq (n+1)$, we see that

$$-K_{\mathbb{P}(V)}(A_2) = (n+1-c_1)a + r \cdot \xi(A_2) \geq 2(n+1-c_1) + r$$
$$> n + r = \dim(\mathbb{P}(V)) + 1;$$

but this contradicts with $-K_{\mathbb{P}(V)}(A_2) \leq \dim(\mathbb{P}(V)) + 1$. Thus $a = 1$ and $A_2 = (h_{n-1}\xi_{r-1} + bh_n\xi_{r-2})_*$. Now $[\pi(E_2)] = \pi_*(A_2) = (h_{n-1})_*$. So $\pi(E_2)$ is a line in $\mathbb{P}^n$. Since $V|_\ell = \oplus_{i=1}^r \mathcal{O}_\ell(m_i)$ for a generic line $\ell \subset \mathbb{P}^n$, $V|_{\pi(E_2)} = \oplus_{i=1}^r \mathcal{O}_{\pi(E_2)}(m_i')$ where $m_i' \geq m_1$ for every $i$. Thus, $\xi(A_2) \geq m_1$, and so $c_1 + b \geq m_1$. It follows that

$$A_2 = [h_{n-1}\xi_{r-1} + (m_1 - c_1)h_n\xi_{r-2}]_* + (c_1 + b - m_1) \cdot (h_n\xi_{r-2})_*.$$

Therefore, $A_2 = [h_{n-1}\xi_{r-1} + (m_1 - c_1)h_n\xi_{r-2}]_* = A$.

(iv) Since $\xi(A) = m_1 = 1$ and $\xi$ is ample, the conclusion follows. $\square$

Next, let $\mathfrak{M}(A, 0)$ be the moduli space of morphisms $f : \mathbb{P}^1 \to \mathbb{P}(V)$ with $[\text{Im}(f)] = A$. In the lemma below, we study the morphisms in $\mathfrak{M}(A, 0)$ when $A = [h_{n-1}\xi_{r-1} + (m - c_1)h_n\xi_{r-2}]_*$. Note that $\xi(A) = m$.

**Lemma 2.3.** *Let $A = [h_{n-1}\xi_{r-1} + (m - c_1)h_n\xi_{r-2}]_*$.*

(i) *If $\mathfrak{M}(A, 0) \neq \emptyset$, then $m \geq m_1$ and $\mathfrak{M}(A, 0)$ consists of embeddings $f : \ell \to \mathbb{P}(V)$ induced by surjective maps $V|_\ell \to \mathcal{O}_\ell(m) \to 0$ where $\ell$ are lines in $\mathbb{P}^n$;*

(ii) *If $m = m_1$ and $m_1 = \ldots = m_k < m_{k+1} \leq \ldots \leq m_r$, then the moduli space $\mathfrak{M}(A, 0)$ has (complex) dimension $(2n + k)$;*

(iii) *If $m \geq m_r$, then $\mathfrak{M}(A, 0)$ has dimension $(2n + r + rm - c_1)$.*

*Proof.* (i) Let $f : \mathbb{P}^1 \to \mathbb{P}(V)$ be a morphism in $\mathfrak{M}(A, 0)$. Then $[\text{Im}(f)] = A = [h_{n-1}\xi_{r-1} + (m - c_1)h_n\xi_{r-2}]_*$. Since $h(A) = 1$, $\pi^*H \cap f(\mathbb{P}^1)$ consists of a single point for any hyperplane $H$ in $\mathbb{P}^n$. Thus, $\pi|_{f(\mathbb{P}^1)} : f(\mathbb{P}^1) \to (\pi \circ f)(\mathbb{P}^1)$ is an isomorphism and $\ell = (\pi \circ f)(\mathbb{P}^1)$ is a line in $\mathbb{P}^n$. Since $h([\ell]) = 1$, $(\pi \circ f) : \mathbb{P}^1 \to \ell = (\pi \circ f)(\mathbb{P}^1)$ is also an isomorphism, and so is $f : \mathbb{P}^1 \to f(\mathbb{P}^1)$. Replacing $f : \mathbb{P}^1 \to \mathbb{P}(V)$ by $f \circ (\pi \circ f)^{-1} : \ell \to \mathbb{P}(V)$, we conclude that $\mathfrak{M}(A, 0)$ consists of embeddings $f : \ell \to \mathbb{P}(V)$ such that $[\text{Im}(f)] = A$, $\ell$ are lines in $\mathbb{P}^n$, and $\pi|_{f(\ell)} : f(\ell) \to \ell$ are isomorphisms. In particular, these embeddings $f : \ell \to \mathbb{P}(V)$ are sections to the natural projection $\pi|_{\mathbb{P}(V|_\ell)} : \mathbb{P}(V|_\ell) \to \ell$. Thus, by the Proposition 7.12 in Chapter II of [7], these embeddings are induced by surjective maps $V|_\ell \to \mathcal{O}_\ell(m) \to 0$. By (2.1), the splitting type of the restrictions of $V$ to generic lines in $\mathbb{P}^n$ is $(m_1, \ldots, m_r)$ with $m_1 \leq \ldots \leq m_r$; thus we must have $V|_\ell = \oplus_{i=1}^r \mathcal{O}_\ell(m_i')$ where $m_i' \geq m_1$ for every $i$. It follows that $m \geq \min\{m_1', \ldots, m_r'\} \geq m_1$.

(ii) Note that all the lines in $\mathbb{P}^n$ are parameterized by the Grassmannian $G(2, n+1)$ which has dimension $2(n-1)$. For a fixed generic line $\ell \subset \mathbb{P}^n$, the surjective maps $V|_\ell \to \mathcal{O}_\ell(m_1) \to 0$ are parameterized by

$$\mathbb{P}(\text{Hom}(V|_\ell, \mathcal{O}_\ell(m_1))) \cong \mathbb{P}(\oplus_{i=1}^r H^0(\ell, \mathcal{O}_\ell(m_1 - m_i))) \cong \mathbb{P}^{k-1};$$

It follows from (i) that as the generic line $\ell$ varies, the morphisms $f : \ell \to \mathbb{P}(V)$ induced by these surjective maps $V|_\ell \to \mathcal{O}_\ell(m_1) \to 0$ form an open dense subset of $\mathfrak{M}(A, 0)$. Thus, $\dim(\mathfrak{M}(A, 0)) = 3 + 2(n - 1) + (k - 1) = 2n + k$.

(iii) As in the proof of (ii), for a fixed generic line $\ell \subset \mathbb{P}^n$, the surjective maps $V|_\ell \to \mathcal{O}_\ell(m) \to 0$ are parameterized by a nonempty open subset of

$$\mathbb{P}(\text{Hom}(V|_\ell, \mathcal{O}_\ell(m))) \cong \mathbb{P}(\oplus_{i=1}^r H^0(\ell, \mathcal{O}_\ell(m - m_i))) \cong \mathbb{P}^{(rm-c_1+r)-1}.$$

As the generic line $\ell$ varies, the morphisms $f : \ell \to \mathbb{P}(V)$ induced by these surjective maps $V|_\ell \to \mathcal{O}_\ell(m) \to 0$ form an open dense subset of $\mathfrak{M}(A, 0)$. It follows that $\mathfrak{M}(A, 0)$ has dimension $(2n + r + rm - c_1)$. $\square$

## 3. Calculation of Gromov-Witten invariants

In this section, we shall compute some Gromov-Witten invariants of $\mathbb{P}(V)$. First of all, we recall that for two homogeneous elements $\alpha$ and $\beta$ in $H^*(\mathbb{P}(V); \mathbb{Z})$, the

quantum product $\alpha \cdot \beta \in H^*(\mathbb{P}(V); \mathbb{Z})$ can be written as

$$\alpha \cdot \beta = \sum_{A \in H_2(\mathbb{P}(V); \mathbb{Z})} (\alpha \cdot \beta)_A \cdot e^{t \cdot K_{\mathbb{P}(V)}(A)} \tag{3.1}$$

where $(\alpha \cdot \beta)_A$ has degree $\deg(\alpha) + \deg(\beta) + 2K_{\mathbb{P}(V)}(A)$ and is defined by

$$(\alpha \cdot \beta)_A(\gamma_*) = \Phi_{(A,0)}(\alpha, \beta, \gamma)$$

for a homogeneous cohomology class $\gamma \in H^*(\mathbb{P}(V); \mathbb{Z})$ with

$$\deg(\gamma) = -2K_{\mathbb{P}(V)}(A) + 2(n + r - 1) - \deg(\alpha) - \deg(\beta). \tag{3.2}$$

Furthermore, for higher quantum products, we have

$$\alpha_1 \cdot \alpha_2 \cdot \ldots \cdot \alpha_k = \sum_{A \in H_2(\mathbb{P}(V); \mathbb{Z})} (\alpha_1 \cdot \alpha_2 \cdot \ldots \cdot \alpha_k)_A \cdot e^{t \cdot K_{\mathbb{P}(V)}(A)} \tag{3.3}$$

where $(\alpha_1 \cdot \alpha_2 \cdot \ldots \cdot \alpha_k)_A$ is defined as $(\alpha_1 \cdot \alpha_2 \cdot \ldots \cdot \alpha_k)_A(\gamma_*) = \Phi_{(A,0)}(\alpha_1, \alpha_2, \ldots, \alpha_k, \gamma)$. Thus, $\alpha_1 \cdot \alpha_2 \cdot \ldots \cdot \alpha_k = \alpha_1 \alpha_2 \ldots \alpha_k +$ (lower order terms), where $\alpha_1 \alpha_2 \ldots \alpha_k$ stands for the ordinary cohomology product of $\alpha_1, \alpha_2, \ldots, \alpha_k$, and the degree of a lower order term is dropped by $2K_{\mathbb{P}(V)}(A)$ for some $A \in H_2(\mathbb{P}(V); \mathbb{Z})$ which is represented by a nonconstant effective rational curve.

There are two explanations for the Gromov-Witten invariant $\Phi_{(A,0)}(\alpha, \beta, \gamma)$ defined by the second author [12]. Recall that the Gromov-Witten invariant is only defined for a generic almost complex structure and that $\mathfrak{M}(A, 0)$ is the moduli space of morphisms $f : \mathbb{P}^1 \to \mathbb{P}(V)$ with $[\text{Im}(f)] = A$. Assume the genericity conditions:

(i) $\mathfrak{M}(A, 0)/PSL(2; \mathbb{C})$ is smooth in the sense that $h^1(N_f) = 0$ for every $f \in \mathfrak{M}(A, 0)$ where $N_f$ is the normal bundle, and

(ii) the homology class $A$ is only represented by irreducible and reduced curves.

Then the complex structure is already generic and one can use algebraic geometry to calculate the Gromov-Witten invariants. Moreover, $\mathfrak{M}(A, 0)/PSL(2; \mathbb{C})$ is compact with the expected complex dimension

$$-K_{\mathbb{P}(V)}(A) + (n + r - 1) - 3. \tag{3.4}$$

The first explanation for $\Phi_{(A,0)}(\alpha, \beta, \gamma)$ is that when $\alpha, \beta, \gamma$ are classes of subvarieties $B, C, D$ of $\mathbb{P}(V)$ in general position, $\Phi_{(A,0)}(\alpha, \beta, \gamma)$ is the number of rational curves $E$ in $\mathbb{P}(V)$ such that $[E] = A$ and that $E$ intersects with $B, C, D$ (counted with suitable multiplicity). The second explanation for $\Phi_{(A,0)}(\alpha, \beta, \gamma)$ is that

$$\Phi_{(A,0)}(\alpha, \beta, \gamma) = \int_{\mathfrak{M}(A,0)} e_0^*(\alpha) \cdot e_1^*(\beta) \cdot e_2^*(\gamma)$$

where the evaluation map $e_i : \mathfrak{M}(A, 0) \to \mathbb{P}(V)$ is defined by $e_i(f) = f(i)$.

Assume that the genericity condition (i) is not satisfied but $h^1(N_f)$ is independent of $f \in \mathfrak{M}(A, 0)$ and $\mathfrak{M}(A, 0)/PSL(2; \mathbb{C})$ is smooth with dimension

$$-K_{\mathbb{P}(V)}(A) + (n + r - 1) - 3 + h^1(N_f).$$

Then one can form an obstruction bundle $COB$ of rank $h^1(N_f)$ over the moduli space $\mathfrak{M}(A,0)$. Moreover, if the genericity condition (ii) is satisfied, then by the Proposition 5.7 in [11], we have

$$\Phi_{(A,0)}(\alpha,\beta,\gamma) = \int_{\mathfrak{M}(A,0)} e_0^*(\alpha) \cdot e_1^*(\beta) \cdot e_2^*(\gamma) \cdot e(COB) \qquad (3.5)$$

where $e(COB)$ stands for the Euler class of the bundle $COB$.

We remark that in general, the cohomology class $h_i\xi_j$ may not be able to be represented by a subvariety of $\mathbb{P}(V)$. However, since $\xi$ is ample, $s\xi$ is very ample for $s \gg 0$. Thus, the multiple $th_i\xi_j$ with $t \gg 0$ can be represented by a subvariety of $\mathbb{P}(V)$ whose image in $\mathbb{P}^n$ is a linear subspace of codimension $i$. Since $\Phi_{(A,0)}(\alpha,\beta,h_i\xi_j) = 1/t \cdot \Phi_{(A,0)}(\alpha,\beta,t \cdot h_i\xi_j)$ for $\alpha$ and $\beta$ in $H^*(\mathbb{P}(V);\mathbb{Z})$, it follows that to compute $\Phi_{(A,0)}(\alpha,\beta,h_i\xi_j)$, it suffices to compute $\Phi_{(A,0)}(\alpha,\beta,t \cdot h_i\xi_j)$. In the proofs below, we shall assume implicitly that $t = 1$ for simplicity.

Now we compute the Gromov-Witten invariant $\Phi_{((h_n\xi_{r-2})_*,0)}(\xi,\xi_{r-1},h_n\xi_{r-1})$.

**Lemma 3.6.** $\Phi_{((h_n\xi_{r-2})_*,0)}(\xi,\xi_{r-1},h_n\xi_{r-1}) = 1$.

*Proof.* First of all, we notice that $A = (h_n\xi_{r-2})_*$ can only be represented by lines $\ell$ in the fibers of $\pi$. In particular, there is no reducible or nonreduced effective curves representing $A$. Thus, $\mathfrak{M}(A,0)/PSL(2;\mathbb{C})$ is compact and has dimension:

$$\dim(\mathbb{P}^n) + \dim G(2,r) = n + 2(r-2) = n + 2r - 4$$

which is the expected dimension by (3.4) (here we use $G(2,r)$ to stand for the Grassmannian of lines in $\mathbb{P}^{r-1}$). Next, we want to show that $\mathfrak{M}(A,0)/PSL(2;\mathbb{C})$ is smooth. Let $p = \pi(\ell)$. Then from the two inclusions $\ell \subset \pi^{-1}(p) \subset \mathbb{P}(V)$, we obtain an exact sequence relating normal bundles:

$$0 \to N_{\ell|\pi^{-1}(p)} \to N_{\ell|\mathbb{P}(V)} \to (N_{\pi^{-1}(p)|\mathbb{P}(V)})|_\ell \to 0.$$

Since $N_{\ell|\pi^{-1}(p)} = N_{\ell|\mathbb{P}^{r-1}} = \mathcal{O}_\ell(1)^{\oplus(r-2)}$ and $N_{\pi^{-1}(p)|\mathbb{P}(V)} = (\pi|_{\pi^{-1}(p)})^*T_{p,\mathbb{P}^n}$, the previous exact sequence is simplified into the exact sequence

$$0 \to \mathcal{O}_\ell(1)^{\oplus(r-2)} \to N_{\ell|\mathbb{P}(V)} \to (\pi|_\ell)^*T_{p,\mathbb{P}^n} \to 0.$$

It follows that $H^1(\ell, N_{\ell|\mathbb{P}(V)}) = 0$. Thus, $\mathfrak{M}(A,0)/PSL(2;\mathbb{C})$ is smooth.

Finally, the Poincaré dual of $h_n\xi_{r-1}$ is represented by a point $q_0 \in \mathbb{P}(V)$. If a line $\ell \in \mathfrak{M}(A,0)$ intersects $q_0$, then $\ell \subset \pi^{-1}(\pi(q_0))$. Since the restriction of $\xi$ to the fiber $\pi^{-1}(\pi(q_0)) \cong \mathbb{P}^{r-1}$ is the cohomology class of a hyperplane in $\mathbb{P}^{r-1}$, we conclude that $\Phi_{((h_n\xi_{r-2})_*,0)}(\xi,\xi_{r-1},h_n\xi_{r-1}) = 1$. $\square$

Next we show the vanishing of some Gromov-Witten invariant.

**Lemma 3.7.** Let $A = b(h_n\xi_{r-2})_*$ with $b \geq 1$ and $\alpha \in H^*(\mathbb{P}(V);\mathbb{Z})$. Then,

$$\Phi_{(A,0)}(h_{p_1}\xi_{q_1}, h_{p_2}\xi_{q_2}, \alpha) = 0$$

if $p_1, q_1, p_2, q_2$ are nonnegative integers with $(q_1 + q_2) < r$.

*Proof.* We may assume that $\alpha$ is a homogeneous class in $H^*(\mathbb{P}(V); \mathbb{Z})$. By (3.2),

$$\frac{1}{2} \cdot \deg(\alpha) = (n + r - 1) - K_{\mathbb{P}(V)}(A) - (p_1 + p_2 + q_1 + q_2)$$
$$= (n + r + br - 1) - (p_1 + p_2 + q_1 + q_2).$$

Let $\alpha = h_{(n+r+br-1)-(p_1+p_2+q_1+q_2+q_3)}\xi_{q_3}$ with $0 \leq q_3 \leq (r-1)$. Let $B, C, D$ be the subvarieties of $\mathbb{P}(V)$ in general position, whose homology classes are Poincaré dual to $h_{p_1}\xi_{q_1}, h_{p_2}\xi_{q_2}, \alpha$ respectively. Then the homology classes of $\pi(B), \pi(C), \pi(D)$ in $\mathbb{P}^n$ are Poincaré dual to $h_{p_1}, h_{p_2}, h_{(n+r+br-1)-(p_1+p_2+q_1+q_2+q_3)}$ respectively. Since $(q_1 + q_2 + q_3) < (2r-1)$, we have $p_1 + p_2 + [(n+r+br-1) - (p_1+p_2+q_1+q_2+q_3)] = (n+r+br-1) - (q_1+q_2+q_3) > n$. Thus, $\pi(B) \cap \pi(C) \cap \pi(D) = \emptyset$. Notice that the genericity conditions (i) and (ii) mentioned earlier in this section are not satisfied for $b \geq 2$. However, we observe that these conditions can be relaxed by assuming:

(i') $h^1(N_f) = 0$ for every $f \in \mathfrak{M}(A, 0)$ such that $\text{Im}(f)$ intersects $B, C, D$, and
(ii') there is no reducible or nonreduced effective (connected) curve $E$ such that $[E] = A$ and $E$ intersects $B, C, D$.

In fact, we will show that there is no effective connected curve $E$ at all representing $A$ and intersecting $B, C, D$. It obviously implies (i'), (ii') and

$$\Phi_{(A,0)}(h_{p_1}\xi_{q_1}, h_{p_2}\xi_{q_2}, \alpha) = 0.$$

Suppose that $E = \sum a_i E_i$ is such an effective connected curve where $a_i > 0$ and $E_i$ is irreducible and reduced. Then, $\sum a_i[E_i] = [E] = A$. Since $(h_n\xi_{r-2})_*$ generates an extremal ray for $\mathbb{P}(V)$, $[E_i] = b_i(h_n\xi_{r-2})_*$ for $0 < b_i \leq b$. Thus the curves $E_i$ are contained in the fibers of $\pi$. Since $E$ is connected, all the curves $E_i$ must be contained in the same fiber of $\pi$. So $\pi(E)$ is a single point. Since $E$ intersects $B, C, D$, $\pi(E)$ intersects with $\pi(B), \pi(C), \pi(D)$. It follows that $\pi(B) \cap \pi(C) \cap \pi(D)$ contains $\pi(E)$ and is nonempty. Therefore we obtain a contradiction. $\square$

Finally, we show that if $c_1 < 2r$ and $A = [h_{n-1}\xi_{r-1} + (1-c_1)h_n\xi_{r-2}]_*$, then $\Phi_{(A,0)}(h, h_n, h_n\xi_{2r-c_1-1}) = 1$. Since $c_1 < 2r$, we see that for a generic line $\ell \subset \mathbb{P}^n$,

$$V|_\ell = \mathcal{O}_\ell(1)^{\oplus k} \oplus \mathcal{O}_\ell(m_{k+1}) \oplus \ldots \oplus \mathcal{O}_\ell(m_r)$$

where $k \geq 1$ and $2 \leq m_{k+1} \leq \ldots \leq m_r$. We remark that even though the moduli space $\mathfrak{M}(A, 0)/PSL(2; \mathbb{C})$ is compact by Lemma 2.2 (iv), it may not have the correct dimension by Lemma 2.3 (ii). The proof is lengthy, but the basic idea is that we shall determine the obstruction bundle and use the formula (3.5).

**Lemma 3.8.** *Let $V$ be a rank-$r$ ample vector bundle over $\mathbb{P}^n$ satisfying $c_1 < 2r$ and the assumption of Theorem B (i). If $A = [h_{n-1}\xi_{r-1} + (1-c_1)h_n\xi_{r-2}]_*$, then*

$$\Phi_{(A,0)}(h, h_n, h_n\xi_{2r-c_1-1}) = 1.$$

*Proof.* Note that by Lemma 2.2 (iv), the moduli space $\mathfrak{M}(A, 0)/PSL(2; \mathbb{C})$ is compact. Let $B, C, D$ be the subvarieties of $\mathbb{P}(V)$ in general position, whose homology classes are Poincaré dual to $h, h_n, h_n \xi_{2r-c_1-1}$ respectively. Then the homology classes of $\pi(B), \pi(C), \pi(D)$ in $\mathbb{P}^n$ are Poincaré dual to $h, h_n, h_n$ respectively. Thus $\pi(C)$ and $\pi(D)$ are two different points in $\mathbb{P}^n$. Let $\ell_0$ be the unique line passing $\pi(C)$ and $\pi(D)$. Let $V|_{\ell_0} = \mathcal{O}_{\ell_0}(1)^{\oplus k} \oplus \mathcal{O}_{\ell_0}(m_{k+1}) \oplus \ldots \oplus \mathcal{O}_{\ell_0}(m_r)$ where $2 \leq m_{k+1} \leq \ldots \leq m_r$. Since $c_1 < 2r$, $k \geq 1$. Let $f : \ell \to \mathbb{P}(V)$ be a morphism in $\mathfrak{M}(A,0)$ for some line $\ell \in \mathbb{P}^n$. If $\text{Im}(f)$ intersects with $B, C$, and $D$, then $\ell = \ell_0$. As in the proof of Lemma 2.3 (ii), the morphisms $f : \ell_0 \to \mathbb{P}(V)$ in $\mathfrak{M}(A,0)$ are parameterized by $\mathbb{P}(\text{Hom}(V|_{\ell_0}, \mathcal{O}_{\ell_0}(1))) \cong \mathbb{P}^{k-1}$; moreover, $\text{Im}(f)$ are of the form:

$$\ell_0 \times \{q\} \subset \ell_0 \times \mathbb{P}^{k-1} = \mathbb{P}(\mathcal{O}_{\ell_0}(1)^{\oplus k}) \subset \mathbb{P}(V|_{\ell_0}) \subset \mathbb{P}(V) \qquad (3.9)$$

where $q$ stands for points in $\mathbb{P}^{k-1} \subset \mathbb{P}^{r-1} \cong \pi^{-1}(\pi(D))$. Note that $\ell_0 \times \{q\}$ always intersects with $B$ and $C$, and that $D$ is a dimension-$(c_1 - r)$ linear subspace in $\mathbb{P}^{r-1} \cong \pi^{-1}(\pi(D))$. Thus, $\ell_0 \times \{q\}$ intersects with $B, C, D$ simultaneously if and only if $\ell_0 \times \{q\}$ intersects with $D$, and if only only if

$$q \in \mathbb{P}^{c_1+k-2r} \stackrel{\text{def}}{=} \mathbb{P}^{k-1} \cap D \subset \mathbb{P}^{r-1} \cong \pi^{-1}(\pi(D)). \qquad (3.10)$$

It follows that $\mathfrak{M}/PSL(2;\mathbb{C}) \cong \mathbb{P}^{c_1+k-2r}$ where $\mathfrak{M}$ consists of morphisms $f \in \mathfrak{M}(A,0)$ such that $\text{Im}(f)$ intersects with $B, C, D$ simultaneously.

If $c_1 + k - 2r = 0$, then $a_0 = \Phi_{(A,0)}(h, h_n, h_n\xi_{2r-c_1-1}) = 1$. But in general, we have $c_1 + k - 2r \geq 0$. We shall use (3.5) to compute $a_0 = \Phi_{(A,0)}(h, h_n, h_n\xi_{2r-c_1-1})$. Let $N_f = N_{\ell_0 \times \{q\}|\mathbb{P}(V)}$ be the normal bundle of $\text{Im}(f) = \ell_0 \times \{q\}$ in $\mathbb{P}(V)$. If $h^1(N_f)$ is constant for every $f \in \mathfrak{M}$, then by (3.5), $\Phi_{(A,0)}(h, h_n, h_n\xi_{2r-c_1-1})$ is the Euler number $e(COB)$ of the rank-$(c_1 + k - 2r)$ obstruction bundle $COB$ over

$$\mathfrak{M}/PSL(2;\mathbb{C}) \cong \mathbb{P}^{c_1+k-2r}.$$

Thus we need to show that $h^1(N_f)$ is constant for every $f \in \mathfrak{M}$.

First, we study the normal bundle $N_{\ell_0 \times \mathbb{P}^{c_1+k-2r}|\mathbb{P}(V)}$. The three inclusions

$$\ell_0 \times \mathbb{P}^{c_1+k-2r} \subset \ell_0 \times \mathbb{P}^{k-1} = \mathbb{P}(\mathcal{O}_{\ell_0}(1)^{\oplus k}) \subset \mathbb{P}(V|_{\ell_0}) \subset \mathbb{P}(V) \qquad (3.11)$$

give rise to two exact sequences relating normal bundles:

$$0 \to N_{\ell_0 \times \mathbb{P}^{c_1+k-2r}|\mathbb{P}(V|_{\ell_0})} \to N_{\ell_0 \times \mathbb{P}^{c_1+k-2r}|\mathbb{P}(V)} \to N_{\mathbb{P}(V|_{\ell_0})|\mathbb{P}(V)} \to 0$$

$$0 \to N_{\ell_0 \times \mathbb{P}^{c_1+k-2r}|\mathbb{P}(\mathcal{O}_{\ell_0}(1)^{\oplus k})} \to N_{\ell_0 \times \mathbb{P}^{c_1+k-2r}|\mathbb{P}(V|_{\ell_0})} \to N_{\mathbb{P}(\mathcal{O}_{\ell_0}(1)^{\oplus k})|\mathbb{P}(V|_{\ell_0})} \to 0$$

Notice that $N_{\mathbb{P}(V|_{\ell_0})|\mathbb{P}(V)} = (\pi|_{\mathbb{P}(V|_{\ell_0})})^*(N_{\ell_0|\mathbb{P}^n}) = \mathcal{O}_{\ell_0}(1)^{\oplus(n-1)}$ and that

$$N_{\ell_0 \times \mathbb{P}^{c_1+k-2r}|\mathbb{P}(\mathcal{O}_{\ell_0}(1)^{\oplus k})} = N_{\ell_0 \times \mathbb{P}^{c_1+k-2r}|\ell_0 \times \mathbb{P}^{k-1}} = \mathcal{O}_{\mathbb{P}^{c_1+k-2r}}(1)^{\oplus(2r-c_1-1)}.$$

Since $V|_{\ell_0} = \mathcal{O}_{\ell_0}(1)^{\oplus k} \oplus \oplus_{i=k+1}^{r} \mathcal{O}_{\ell_0}(m_i)$, $\xi|_{\ell_0 \times \mathbb{P}^{k-1}} = \mathcal{O}_{\ell_0}(1) \otimes \mathcal{O}_{\mathbb{P}^{k-1}}(1)$ and

$$N_{\mathbb{P}(\mathcal{O}_{\ell_0}(1)^{\oplus k})|\mathbb{P}(V|_{\ell_0})} = \oplus_{i=k+1}^{r} \mathcal{O}_{\ell_0}(-m_i) \otimes \xi|_{\ell_0 \times \mathbb{P}^{k-1}}$$
$$= \oplus_{i=k+1}^{r} \mathcal{O}_{\ell_0}(1-m_i) \otimes \mathcal{O}_{\mathbb{P}^{k-1}}(1).$$

Thus the previous two exact sequences are simplified to:

$$0 \to N_{\ell_0 \times \mathbb{P}^{c_1+k-2r}|\mathbb{P}(V|_{\ell_0})} \to N_{\ell_0 \times \mathbb{P}^{c_1+k-2r}|\mathbb{P}(V)} \to \mathcal{O}_{\ell_0}(1)^{\oplus(n-1)} \to 0 \qquad (3.12)$$

$$0 \to \mathcal{O}_{\mathbb{P}^{c_1+k-2r}}(1)^{\oplus(2r-c_1-1)} \to N_{\ell_0 \times \mathbb{P}^{c_1+k-2r}|\mathbb{P}(V|_{\ell_0})} \to$$
$$\oplus_{i=k+1}^{r} \mathcal{O}_{\ell_0}(1-m_i) \otimes \mathcal{O}_{\mathbb{P}^{c_1+k-2r}}(1) \to 0 \qquad (3.13)$$

Now (3.13) splits since for $k+1 \leq i \leq r$, we have $m_i \geq 2$ and

$$\mathrm{Ext}^1(\mathcal{O}_{\ell_0}(1-m_i) \otimes \mathcal{O}_{\mathbb{P}^{c_1+k-2r}}(1), \mathcal{O}_{\mathbb{P}^{c_1+k-2r}}(1))$$
$$= H^1(\ell_0 \times \mathbb{P}^{c_1+k-2r}, \mathcal{O}_{\ell_0}(m_i - 1)) = 0.$$

Thus, the normal bundle $N_{\ell_0 \times \mathbb{P}^{c_1+k-2r}|\mathbb{P}(V|_{\ell_0})}$ is isomorphic to

$$\oplus_{i=k+1}^{r} \mathcal{O}_{\ell_0}(1-m_i) \otimes \mathcal{O}_{\mathbb{P}^{c_1+k-2r}}(1) \oplus \mathcal{O}_{\mathbb{P}^{c_1+k-2r}}(1)^{\oplus(2r-c_1-1)},$$

and the exact sequence (3.12) becomes to the exact sequence:

$$0 \to \oplus_{i=k+1}^{r} \mathcal{O}_{\ell_0}(1-m_i) \otimes \mathcal{O}_{\mathbb{P}^{c_1+k-2r}}(1) \oplus \mathcal{O}_{\mathbb{P}^{c_1+k-2r}}(1)^{\oplus(2r-c_1-1)} \to$$
$$N_{\ell_0 \times \mathbb{P}^{c_1+k-2r}|\mathbb{P}(V)} \to \mathcal{O}_{\ell_0}(1)^{\oplus(n-1)} \to 0 \qquad (3.14)$$

Restricting (3.14) to $\ell_0 \times \{q\}$ and taking long exact cohomology sequence result

$$\oplus_{i=k+1}^{r} H^1(\mathcal{O}_{\ell_0}(1-m_i)) \otimes \mathcal{O}_{\mathbb{P}^{c_1+k-2r}}(1)|_q \to$$
$$H^1((N_{\ell_0 \times \mathbb{P}^{c_1+k-2r}|\mathbb{P}(V)})|_{\ell_0 \times \{q\}}) \to 0. \qquad (3.15)$$

Next, we determine $N_f$ and show that $h^1(N_f) \leq c_1 + k - 2r$. The two inclusions $\ell_0 \times \{q\} \subset \ell_0 \times \mathbb{P}^{c_1+k-2r} \subset \mathbb{P}(V)$ give an exact sequence

$$0 \to N_{\ell_0 \times \{q\}|\ell_0 \times \mathbb{P}^{c_1+k-2r}} \to N_{\ell_0 \times \{q\}|\mathbb{P}(V)} \to (N_{\ell_0 \times \mathbb{P}^{c_1+k-2r}|\mathbb{P}(V)})|_{\ell_0 \times \{q\}} \to 0.$$

Since $N_{\ell_0 \times \{q\}|\ell_0 \times \mathbb{P}^{c_1+k-2r}} = T_{q,\mathbb{P}^{c_1+k-2r}}$, the above exact sequence becomes

$$0 \to T_{q,\mathbb{P}^{c_1+k-2r}} \to N_f \to (N_{\ell_0 \times \mathbb{P}^{c_1+k-2r}|\mathbb{P}(V)})|_{\ell_0 \times \{q\}} \to 0. \qquad (3.16)$$

Thus, $h^1(N_f) = h^1((N_{\ell_0 \times \mathbb{P}^{c_1+k-2r}|\mathbb{P}(V)})|_{\ell_0 \times \{q\}})$. By (3.15), we obtain

$$h^1(N_f) = h^1((N_{\ell_0 \times \mathbb{P}^{c_1+k-2r}|\mathbb{P}(V)})|_{\ell_0 \times \{q\}}) \leq \sum_{i=k+1}^{r} h^1(\mathcal{O}_{\ell_0}(1-m_i))$$
$$= \sum_{i=k+1}^{r} (m_i - 2) = c_1 + k - 2r.$$

Finally, we show that $h^1(N_f) = c_1 + k - 2r$. It suffices to prove that $h^1(N_f) \geq c_1 + k - 2r$. Since $\ell_0$ is a generic line in $\mathbb{P}^n$ and $V|_{\ell_0} = \mathcal{O}_{\ell_0}(1)^{\oplus k} \oplus \oplus_{i=k+1}^{r} \mathcal{O}_{\ell_0}(m_i)$, $\dim \mathfrak{M}(A, 0) = (2n + k)$ by Lemma 2.3 (ii). Since $h^0(N_f)$ is the dimension of the Zariski tangent space of $\mathfrak{M}(A, 0)/PSL(2;\mathbb{C})$ at $f$, $h^0(N_f) \geq (2n + k - 3)$. Thus,

$$h^1(N_f) = h^0(N_f) - \chi(N_f) \geq (2n + k - 3) - (2n + 2r - c_1 - 3) = k + c_1 - 2r.$$

Therefore, $h^1(N_f) = c_1 + k - 2r$. In particular, $h^1(N_f)$ is independent of $f \in \mathfrak{M}$. To obtain the obstruction bundle $COB$ over $\mathbb{P}^{c_1+k-2r}$, we notice that (3.15) gives

$$\oplus_{i=k+1}^{r} H^1(\mathcal{O}_{\ell_0}(1 - m_i)) \otimes \mathcal{O}_{\mathbb{P}^{c_1+k-2r}}(1)|_q \cong H^1((N_{\ell_0 \times \mathbb{P}^{c_1+k-2r}|\mathbb{P}(V)})|_{\ell_0 \times \{q\}}).$$

Thus by the exact sequence (3.16), we conclude that

$$\begin{aligned} H^1(N_f) &\cong H^1((N_{\ell_0 \times \mathbb{P}^{c_1+k-2r}|\mathbb{P}(V)})|_{\ell_0 \times \{q\}}) \\ &\cong \oplus_{i=k+1}^{r} H^1(\mathcal{O}_{\ell_0}(1 - m_i)) \otimes \mathcal{O}_{\mathbb{P}^{c_1+k-2r}}(1)|_q. \end{aligned} \tag{3.17}$$

It follows that $COB = \mathcal{O}_{\mathbb{P}^{c_1+k-2r}}(1)^{\oplus(c_1+k-2r)}$. By (3.5), we obtain

$$a_0 = \Phi_{(A,0)}(h, h_n, h_n \xi_{2r-c_1-1}) = e(COB) = 1. \quad \square$$

## 4. Proof of Theorem B

In this section, we prove Theorem B which we restate below.

**Theorem 4.1.** (i) *Let $V$ be a rank-$r$ ample bundle over $\mathbb{P}^n$. Assume either $c_1 \leq n$ or $c_1 \leq (n + r)$ and $V \otimes \mathcal{O}_{\mathbb{P}^n}(-1)$ is nef so that $\mathbb{P}(V)$ is Fano. Then the quantum cohomology $H^*_\omega(\mathbb{P}(V); \mathbb{Z})$ is the ring generated by $h$ and $\xi$ with two relations*

$$h^{n+1} = \sum_{i+j \leq (c_1 - r)} a_{i,j} \cdot h^i \cdot \xi^j \cdot e^{-t(n+1-i-j)} \tag{4.2}$$

$$\sum_{i=0}^{r} (-1)^i c_i \cdot h^i \cdot \xi^{r-i} = e^{-tr} + \sum_{i+j \leq (c_1 - n - 1)} b_{i,j} \cdot h^i \cdot \xi^j \cdot e^{-t(r-i-j)} \tag{4.3}$$

*where the coefficients $a_{i,j}$ and $b_{i,j}$ are integers depending on $V$;*

(ii) *If we further assume that $c_1 < 2r$, then the leading coefficient $a_{0, c_1 - r} = 1$.*

*Proof.* (i) First, we determine the first relation $f^1_\omega$ in (1.3). By Lemma 3.7,

$$h \cdot h_p = h_{p+1} + \sum_{A \in H'_2} (h \cdot h_p)_A \cdot e^{tK_{\mathbb{P}(V)}(A)} \tag{4.4}$$

where $p \geq 1$ and $H'_2$ stands for $H_2(\mathbb{P}(V); \mathbb{Z}) - \mathbb{Z} \cdot (h_n \xi_{r-2})_*$. Thus,

$$h^{n-p} \cdot h_{p+1} = h^{n-p+1} \cdot h_p - \sum_{A \in H'_2} h^{n-p} \cdot (h \cdot h_p)_A \cdot e^{tK_{\mathbb{P}(V)}(A)}.$$

If $(h \cdot h_p)_A \neq 0$, then $A = [E]$ for some effective curve $E$. So $a = h(A) \geq 0$. Since $A \in H'_2$, $a \geq 1$. We claim that $-K_{\mathbb{P}(V)}(A) \geq (n + 1 - c_1 + r)$ with equality if and only if $A = [h_{n-1}\xi_{r-1} + (1-c_1)h_n\xi_{r-2}]_* \stackrel{\text{def}}{=} A_2$. Indeed, if $c_1 \leq n$, then $-K_{\mathbb{P}(V)}(A) = (n+1-c_1)a + r \cdot \xi(A) \geq (n+1-c_1+r)$ with equality if and only if $a = \xi(A) = 1$, that is, if and only if $A = A_2$; if $c_1 \leq (n+r)$ and $(\xi - h)$ is nef, then again $-K_{\mathbb{P}(V)}(A) = (n+1+r-c_1)a + r \cdot (\xi - h)(A) \geq (n+1-c_1+r)$ with equality if and only if $a = 1$ and $(\xi - h)(A) = 0$, that is, if and only if $A = A_2$. Thus, $\deg((h \cdot h_p)_A) = 1 + p + K_{\mathbb{P}(V)}(A) \leq (p - n + c_1 - r)$, and $\deg(h^{n-p} \cdot (h \cdot h_p)_A) \leq (c_1 - r)$. Using induction on $p$ and keeping track of the exponential $e^{tK_{\mathbb{P}(V)}(A)}$, we obtain

$$0 = h_{n+1} = h^{n+1} - \sum_{i+j \leq (c_1 - r)} a_{i,j} \cdot h^i \cdot \xi^j \cdot e^{-t(n+1-i-j)}.$$

Therefore, the first relation $f_\omega^1$ for the quantum cohomology ring is:

$$h^{n+1} = \sum_{i+j \leq (c_1 - r)} a_{i,j} \cdot h^i \cdot \xi^j \cdot e^{-t(n+1-i-j)}.$$

Next, we determine the second relation $f_\omega^2$ in (1.3). We need to compute the quantum product $h^i \cdot \xi^{r-i}$ for $0 \leq i \leq r$. First, we calculate the quantum product $\xi^r$. Note that if $A = (bh_n\xi_{r-2})_*$ with $b \geq 1$, then $-K_{\mathbb{P}(V)}(A) = br \geq r$ with $-K_{\mathbb{P}(V)}(A) = r$ if and only if $A = (h_n\xi_{r-2})_* \stackrel{\text{def}}{=} A_1$. Thus for $p \geq 1$,

$$\xi \cdot \xi_p = \begin{cases} \xi_{p+1} + \sum_{A \in H_2'} (\xi \cdot \xi_p)_A \cdot e^{tK_{\mathbb{P}(V)}(A)}, & \text{if } p < r - 1 \\ \xi_r + (\xi \cdot \xi_{r-1})_{A_1} \cdot e^{-tr} + \sum_{A \in H_2'} (\xi \cdot \xi_{r-1})_A \cdot e^{tK_{\mathbb{P}(V)}(A)}, & \text{if } p = r - 1. \end{cases}$$

Note that $(\xi \cdot \xi_{r-1})_{A_1}$ is of degree zero; by Lemma 3.6, we obtain $(\xi \cdot \xi_{r-1})_{A_1} = \Phi_{(A_1,0)}(\xi, \xi_{r-1}, h_n\xi_{r-1}) = 1$. Therefore for $p \geq 1$,

$$\xi \cdot \xi_p = \begin{cases} \xi_{p+1} + \sum_{A \in H_2'} (\xi \cdot \xi_p)_A \cdot e^{tK_{\mathbb{P}(V)}(A)}, & \text{if } p < r - 1 \\ \xi_r + e^{-tr} + \sum_{A \in H_2'} (\xi \cdot \xi_{r-1})_A \cdot e^{tK_{\mathbb{P}(V)}(A)}, & \text{if } p = r - 1. \end{cases} \quad (4.5)$$

Now, for $i \geq 1$ and $j \geq 1$ with $i + j \leq r$, we have

$$h_i \cdot \xi_j = \begin{cases} h_i\xi_j + \sum_{A \in H_2'} (h_i \cdot \xi_j)_A \cdot e^{tK_{\mathbb{P}(V)}(A)}, & \text{if } i + j < r \\ h_i\xi_j + (h_i\xi_j)_{A_1} \cdot e^{-tr} + \sum_{A \in H_2'} (h_i \cdot \xi_j)_A \cdot e^{tK_{\mathbb{P}(V)}(A)}, & \text{if } i + j = r; \end{cases}$$

when $i + j = r$, $(h_i\xi_j)_{A_1}$ is of degree zero; by Lemma 3.7, we have $(h_i \cdot \xi_j)_{A_1} = \Phi_{(A_1,0)}(h_i, \xi_{r-i}, h_n\xi_{r-1}) = 0$. Therefore for $i \geq 1$ and $j \geq 1$ with $i + j \leq r$,

$$h_i \cdot \xi_j = h_i\xi_j + \sum_{A \in H_2'} (h_i \cdot \xi_j)_A \cdot e^{tK_{\mathbb{P}(V)}(A)}. \quad (4.6)$$

From the proof of the first relation $f_\omega^1$, we see that if $\alpha$ and $\beta$ are homogeneous elements in $H^*(\mathbb{P}(V); \mathbb{Z})$ with $\deg(\alpha) + \deg(\beta) = m \leq r$, then $\deg((\alpha \cdot \beta)_A) \leq m - (n + 1 - c_1 + r)$ for $A \in H_2'$. Thus if $\gamma$ is a homogeneous element in $H^*(\mathbb{P}(V); \mathbb{Z})$ with $\deg(\gamma) = r - m$, then $\deg(\gamma \cdot (\xi \cdot \xi_p)_A) \leq (c_1 - n - 1)$. Since $\sum_{i=0}^{r}(-1)^i c_i \cdot h_i \xi_{r-i} = 0$, it follows from (4.4), (4.5), and (4.6) that the second relation $f_\omega^2$ is

$$\sum_{i=0}^{r}(-1)^i c_i \cdot h^i \cdot \xi^{r-i} = e^{-tr} + \sum_{i+j \leq (c_1-n-1)} b_{i,j} \cdot h^i \cdot \xi^j \cdot e^{-t(r-i-j)}.$$

(ii) From the proof of the first relation in (i), we see that $-K_{\mathbb{P}(V)}(A) \geq (n+1-c_1+r)$ with equality if and only if $A = A_2$; moreover, the term $\xi^{c_1-r}$ can only come from the quantum correction $(h \cdot h_n)_{A_2}$. Now

$$(h \cdot h_n)_{A_2} = \left(\sum_{i=0}^{c_1-r} a'_i h_i \xi_{c_1-r-i}\right) \cdot e^{-t(n+1-c_1+r)}$$

where $a'_0 = \Phi_{(A_2,0)}(h, h_n, h_n \xi_{2r-c_1-1})$. Since $c_1 < 2r$, $(c_1 - r) < r$. By (4.4), (4.5), and (4.6), we conclude that $h_i \xi_{c_1-r-i} = h^i \cdot \xi^{c_1-r-i} +$ (lower degree terms). Thus $a_{0,c_1-r} = a'_0 = \Phi_{(A_2,0)}(h, h_n, h_n \xi_{2r-c_1-1})$. By Lemma 3.8, $a_{0,c_1-r} = 1$. □

It is understood that when $c_1 \leq n$, then the summations on the right-hand-sides of the second relations (4.3) and (4.9) below do not exist.

Next, we shall sharpen the results in Theorem 4.1 by imposing additional conditions on $V$. Let $V$ be a rank-$r$ ample vector bundle over $\mathbb{P}^n$. Then $c_1 \geq r$. Thus if $c_1 < 2r$ and if either $2c_1 \leq (n + r)$ or $2c_1 \leq (n + 2r)$ and $V \otimes \mathcal{O}_{\mathbb{P}^n}(-1)$ is nef, then the conditions in Theorem 4.1 are satisfied.

**Corollary 4.7.** (i) *Let $V$ be a rank-$r$ ample vector bundle over $\mathbb{P}^n$ with $c_1 < 2r$. Assume that either $2c_1 \leq (n+r)$ or $2c_1 \leq (n+2r)$ and $V \otimes \mathcal{O}_{\mathbb{P}^n}(-1)$ is nef so that $\mathbb{P}(V)$ is a Fano variety. Then the first relation (4.2) is*

$$h^{n+1} = \left(\sum_{i=0}^{c_1-r} a_i \cdot h^i \cdot \xi^{c_1-r-i}\right) \cdot e^{-t(n+1+r-c_1)} \qquad (4.8)$$

*where the integers $a_i$ depend on $V$. Moreover, $a_0 = 1$.*

(ii) *Let $V$ be a rank-$r$ ample vector bundle over $\mathbb{P}^n$. Assume that $2c_1 \leq (2n+r+1)$ and $V \otimes \mathcal{O}_{\mathbb{P}^n}(-1)$ is nef so that $\mathbb{P}(V)$ is Fano. Then the second relation (4.3) is*

$$\sum_{i=0}^{r}(-1)^i c_i \cdot h^i \cdot \xi^{r-i} = e^{-tr} + \sum_{i=0}^{c_1-n-1} b_i \cdot h^i \cdot \xi^{c_1-n-1-i} \cdot e^{-t(n+1+r-c_1)} \qquad (4.9)$$

*where the integers $b_i$ depend on $V$.*

*Proof.* (i) From the proof of Theorem 4.1 (i), we notice that it suffices to show that the only homology class $A \in H'_2 = H_2(\mathbb{P}(V); \mathbb{Z}) - \mathbb{Z} \cdot (h_n \xi_{r-2})_*$ which has nonzero contributions to the quantum corrections in (4.4) is $A = [h_{n-1}\xi_{r-1} + (1-c_1)h_n\xi_{r-2}]_* \stackrel{\text{def}}{=} A_2$. In other words, if $A = (ah_{n-1}\xi_{r-1} + bh_n\xi_{r-2})_*$ with $a \neq 0$ and if $\Phi_{(A,0)}(h, h_p, \alpha) \neq 0$ for $1 \leq p \leq n$ and $\alpha \in H^*(\mathbb{P}(V); \mathbb{Z})$, then $A = A_2$. First of all, we show that $a = 1$. Suppose $a \neq 1$. Then $a \geq 2$. By (3.2),

$$\frac{1}{2} \cdot \deg(\alpha) = (n + r - 1) - K_{\mathbb{P}(V)}(A) - 1 - p$$
$$= (n + r - 1) + [(n + 1 - c_1)a + r \cdot \xi(A)] - 1 - p$$
$$\geq \dim(\mathbb{P}(V)) + [(n + 1 - c_1)a + r \cdot \xi(A)] - 1 - n.$$

If $2c_1 \leq (n+r)$, then $c_1 \leq n$, and $[(n+1-c_1)a + r \cdot \xi(A)] - 1 - n \geq 2(n+1-c_1) + r - 1 - n > 0$. If $2c_1 \leq (n+2r)$ and $(\xi - h)$ is nef, then $c_1 \leq n+r$, and $[(n+1-c_1)a + r \cdot \xi(A)] - 1 - n = [(n+1+r-c_1)a + r \cdot (\xi-h)(A)] - 1 - n \geq 2(n+1+r-c_1) - 1 - n > 0$. Thus, $[(n+1-c_1)a + r \cdot \xi(A)] - 1 - n > 0$, and so $\deg(\alpha)/2 > \dim(\mathbb{P}(V))$. But this is absurd. Next, we prove that $b = (1-c_1)$, or equivalently, $\xi(A) = 1$. Suppose $\xi(A) \neq 1$. Then $\xi(A) \geq 2$. By (3.2),

$$\frac{1}{2} \cdot \deg(\alpha) = (n+r-1) + [(n+1-c_1) + r \cdot \xi(A)] - 1 - p$$
$$\geq \dim(\mathbb{P}(V)) + [(n+1-c_1) + 2r] - 1 - n$$
$$> \dim(\mathbb{P}(V))$$

since $c_1 < 2r$. But once again this is absurd.

(ii) We follow the previous arguments for (i). Again it suffices to show that if $A = (ah_{n-1}\xi_{r-1} + bh_n\xi_{r-2})_*$ with $a \neq 0$ and if $\Phi_{(A,0)}(\alpha_1, \alpha_2, \alpha) \neq 0$ for some $\alpha_1, \alpha_2, \alpha \in H^*(\mathbb{P}(V); \mathbb{Z})$ with $\deg(\alpha_1) + \deg(\alpha_2) \leq r$, then $A = A_2$. Indeed, if $a \neq 1$ or if $a = 1$ but $\xi(A) \neq 1$, then we must have $\deg(\alpha)/2 > \dim(\mathbb{P}(V))$. But this is impossible. Therefore, $a = 1$ and $\xi(A) = 1$. So $A = A_2$. $\square$

Now we discuss the relation between the quantum corrections and the extremal rays of the Fano variety $\mathbb{P}(V)$. Let $V$ be a rank-$r$ ample vector bundle over $\mathbb{P}^n$ with $c_1 < 2r$ and $2c_1 \leq (n+r)$. By (4.8) and (4.3), the quantum cohomology ring $H^*_\omega(\mathbb{P}(V); \mathbb{Z})$ is the ring generated by $h$ and $\xi$ with two relations

$$h^{n+1} = \left(\sum_{i=0}^{c_1-r} a_i \cdot h^i \cdot \xi^{c_1-r-i}\right) \cdot e^{-t(n+1+r-c_1)} \tag{4.10}$$

$$\sum_{i=0}^r (-1)^i c_i \cdot h^i \cdot \xi^{r-i} = e^{-tr}. \tag{4.11}$$

From the proof of Theorem 4.1 (i), we notice that the quantum correction to the second relation (4.11) comes from the homology class $A_1 = (h_n \xi_{r-2})_*$ which is represented by the lines in the fibers of $\pi : \mathbb{P}(V) \to \mathbb{P}^n$. Also, we notice from the proof of Corollary 4.7 (i) that the quantum correction to the first relation (4.10) comes from the homology class $A_2 = [h_{n-1}\xi_{r-1} + (1-c_1)h_n\xi_{r-2}]_*$; from the proof of Lemma 3.8, $A_2$ can be represented by a smooth rational curve isomorphic to lines in $\mathbb{P}^n$ via $\pi$. Now $A_1$ generates one of the two extremal rays of $\mathbb{P}(V)$. It is unclear whether $A_2$ generates the other extremal ray. By Lemma 2.2 (iii), if we further assume that $2c_1 \leq (n+1)$, then indeed $A_2$ generates the other extremal ray of $\mathbb{P}(V)$. By Lemma 2.2 (ii), $A_2$ generates the other extremal ray of $\mathbb{P}(V)$ if and only if $(\xi - h)$ is nef, that is, $V \otimes \mathcal{O}_{\mathbb{P}^n}(-1)$ is a nef vector bundle over $\mathbb{P}^n$.

## 5. Direct sum of line bundles over $\mathbb{P}^n$

In this section, we partially verify Batyrev's conjecture on the quantum cohomology of projective bundles associated to direct sum of line bundles over $\mathbb{P}^n$. We shall use (3.5) to compute the necessary Gromov-Witten invariants. Our first step is to recall some standard materials for the Grassmannian $G(2, n+1)$ from [3]. Then we determine certain obstruction bundle and its Euler class. Finally we proceed to determine the first and second relations for the quantum cohomology.

On the Grassmannian $G(2, n+1)$, there exists a tautological exact sequence

$$0 \to S \to (\mathcal{O}_{G(2,n+1)})^{\oplus (n+1)} \to Q \to 0 \tag{5.1}$$

where the sub- and quotient bundles $S$ and $Q$ are of rank $2$ and $(n-1)$ respectively. Let $\alpha$ and $\beta$ be the virtual classes such that $\alpha + \beta = -c_1(S)$ and $\alpha\beta = c_2(S)$. Then

$$\mathrm{cl}(\{\ell \in G(2, n+1) | \ell \cap h_p \neq \emptyset\}) = \frac{\alpha^p - \beta^p}{\alpha - \beta} \tag{5.2}$$

where $\mathrm{cl}(\cdot)$ denotes the fundamental class and $h_p$ stands for a fixed linear subspace of $\mathbb{P}^n$ of codimention $p$. If $P(\alpha, \beta)$ is a symmetric homogeneous polynomial of degree $(2n-2)$ (so that $P(\alpha, \beta)$ can be written as a polynomial of maximal degree in the Chern classes of the bundle $S$), then we have

$$\int_{G(2,n+1)} P(\alpha, \beta) = \left( \text{the coefficient of } \alpha^n \beta^n \text{ in } -\frac{1}{2}(\alpha - \beta)^2 P(\alpha, \beta) \right). \tag{5.3}$$

Let $F_n = \{(x, \ell) \in \mathbb{P}^n \times G(2, n+1) | x \in \ell\}$, and $\pi_1$ and $\pi_2$ are the two natural projections from $F_n$ to $\mathbb{P}^n$ and $G(2, n+1)$ respectively. Then $F_n = \mathbb{P}(S^*)$ where $S^*$ is the dual bundle of $S$, and $(\pi_1^* \mathcal{O}_{\mathbb{P}^n}(1))|_{F_n}$ is the tautological line bundle over $F_n$. Let $\mathrm{Sym}^m(S^*)$ be the $m$-th symmetric product of $S^*$. Then for $m \geq 0$,

$$\pi_{2*}(\pi_1^* \mathcal{O}_{\mathbb{P}^n}(m)|_{F_n}) \cong \mathrm{Sym}^m(S^*). \tag{5.4}$$

By the duality theorem for higher direct image sheaves (see p.253 in [7]),

$$R^1 \pi_{2*}(\pi_1^* \mathcal{O}_{\mathbb{P}^n}(-m)|_{F_n}) \cong (\pi_{2*}(\pi_1^* \mathcal{O}_{\mathbb{P}^n}(m-2)|_{F_n}))^* \otimes (\det S^*)^*$$
$$\cong \mathrm{Sym}^{m-2}(S) \otimes (\det S) \tag{5.5}$$

Now, let $V = \oplus_{i=1}^r \mathcal{O}_{\mathbb{P}^n}(m_i)$ where $1 = m_1 = \ldots = m_k < m_{k+1} \leq \ldots \leq m_r$. Assume that $k \geq 1$ and $\mathbb{P}(V)$ is Fano. Then the two extremal rays of $\mathbb{P}(V)$ are generated by the two classes $A_1 = (h_n \xi_{r-2})_*$ and $A_2 = [h_{n-1} \xi_{r-1} + (1-c_1) h_n \xi_{r-2}]_*$. From the proof of Lemma 2.3 (ii), we see that

$$\mathfrak{M}(A_2, 0)/PSL(2; \mathbb{C}) = G(2, n+1) \times \mathbb{P}^{k-1}. \tag{5.6}$$

Let a morphism $f \in \mathfrak{M}(A_2, 0)$ be induced by some surjective map $V|_\ell \to \mathcal{O}_\ell(1) \to 0$ such that the image $\mathrm{Im}(f)$ of $f$ is of the form

$$\mathrm{Im}(f) = \ell \times \{q\} \subset \ell \times \mathbb{P}^{k-1} \subset \mathbb{P}^n \times \mathbb{P}^{k-1}.$$

Then by arguments similar to the proof of (3.17), we have

$$H^1(N_f) \cong \oplus_{u=k+1}^r H^1(\mathcal{O}_\ell(1 - m_u)) \otimes \mathcal{O}_{\mathbb{P}^{k-1}}(1)|_q. \tag{5.7}$$

It follows that the obstruction bundle $COB$ over $\mathfrak{M}(A_2, 0)/PSL(2; \mathbb{C})$ is

$$COB \cong \oplus_{u=k+1}^r R^1 \pi_{2*}(\pi_1^* \mathcal{O}_{\mathbb{P}^n}(1 - m_u)|_{F_n}) \otimes \mathcal{O}_{\mathbb{P}^{k-1}}(1). \tag{5.8}$$

Since $c_1(S) = -(\alpha + \beta)$ and $c_2(S) = \alpha\beta$, we obtain from (5.5) the following.

**Lemma 5.9.** *The Euler class of the obstruction bundle $COB$ is*

$$e(COB) = \prod_{u=k+1}^{r} \prod_{v=0}^{m_u-3} [(1+v)(-\alpha) + (m_u - 2 - v)(-\beta) + \tilde{h}] \qquad (5.10)$$

*where $\tilde{h}$ stands for the hyperplane class in $\mathbb{P}^{k-1}$.* □

Next assuming $c_1 < 2r$, we shall compute the Gromov-Witten invariant

$$W_i \stackrel{\text{def}}{=} \Phi_{(A_2,0)}(h_{\tilde{n}}, h_{n+1-\tilde{n}}, h_{n-i}\xi_{2r-c_1-1+i}) \qquad (5.11)$$

where $0 \le i \le (c_1 - r)$ and $\tilde{n} = \left[\frac{n+1}{2}\right]$ is the largest integer $\le (n+1)/2$.

**Lemma 5.12.** *Assume $c_1 < \min(2r, (n+1+2r)/2)$ and $0 \le i \le (c_1 - r)$. Then $W_i$ is the coefficient of $t^i$ in the power series expansion of*

$$\prod_{u=1}^{r} (1 - m_u t)^{m_u - 2}.$$

*Proof.* Note that the restriction of $\xi$ to $\mathbb{P}^n \times \mathbb{P}^{k-1} = \mathbb{P}(\mathcal{O}_{\mathbb{P}^n}(1)^{\oplus k})$ is $(h + \tilde{h})$. Thus,

$$h_{n-i}\xi_{2r-c_1-1+i}|_{\mathbb{P}^n \times \mathbb{P}^{k-1}} = \sum_{j=0}^{2r-c_1-1+i} \binom{2r-c_1-1+i}{j} h_{n-i+j} \tilde{h}_{2r-c_1-1+i-j}$$

$$= \sum_{j=0}^{i} \binom{2r-c_1-1+i}{j} h_{n-i+j} \tilde{h}_{2r-c_1-1+i-j}.$$

So by (3.5) (replacing $\mathfrak{M}(A_2, 0)$ by $\mathfrak{M}(A_2, 0)/PSL(2; \mathbb{C})$), (5.2), and Lemma 5.9,

$$W_i = \int_{G(2,n+1) \times \mathbb{P}^{k-1}} \tilde{P}(\alpha, \beta) \qquad (5.13)$$

where $\tilde{P}(\alpha, \beta)$ is the symmetric homogeneous polynomial of degree $(2n-2)+(k-1)$:

$$\tilde{P}(\alpha, \beta) = \frac{\alpha^{\tilde{n}} - \beta^{\tilde{n}}}{\alpha - \beta} \cdot \frac{\alpha^{n+1-\tilde{n}} - \beta^{n+1-\tilde{n}}}{\alpha - \beta}$$

$$\cdot \sum_{j=0}^{i} \binom{2r-c_1-1+i}{j} \frac{\alpha^{n-i+j} - \beta^{n-i+j}}{\alpha - \beta} \cdot \tilde{h}_{2r-c_1-1+i-j}$$

$$\cdot \prod_{u=k+1}^{r} \prod_{v=0}^{m_u-3} [(1+v)(-\alpha) + (m_u - 2 - v)(-\beta) + \tilde{h}]$$

$$= \sum_{j=0}^{i} \binom{2r-c_1-1+i}{j} \frac{\alpha^{n+1} - \alpha^{n+1-\tilde{n}}\beta^{\tilde{n}} - \alpha^{\tilde{n}}\beta^{n+1-\tilde{n}} + \beta^{n+1}}{(\alpha - \beta)^2}$$

$$\cdot \sum_{t=0}^{n-i+j-1} \alpha^t \beta^{n-i+j-1-t} \cdot \tilde{h}_{2r-c_1-1+i-j}$$

$$\cdot \prod_{u=k+1}^{r} \prod_{v=0}^{m_u-3} [(1+v)(-\alpha) + (m_u - 2 - v)(-\beta) + \tilde{h}].$$

By (5.3) and (5.13), we conclude from straightforward manipulations that:

$$W_i = \sum_{j=0}^{i} \binom{2r - c_1 - 1 + i}{j} \cdot (-1)^{i-j}$$

$$\cdot \sum_{j_{k+1}+\ldots+j_r=i-j} \prod_{u=k+1}^{r} \binom{m_u - 2}{j_u}(m_u - 1)^{j_u}$$

$$= \sum_{j=0}^{i} \binom{2r - c_1 - 1 + i}{i - j} \cdot (-1)^{j}$$

$$\cdot \sum_{j_{k+1}+\ldots+j_r=j} \prod_{u=k+1}^{r} \binom{m_u - 2}{j_u}(m_u - 1)^{j_u}.$$

Thus $W_i$ is the coefficient of $t^i$ in the polynomial

$$(1+t)^{2r-c_1-1+i} \cdot \prod_{u=k+1}^{r} [1 - (m_u - 1)t]^{m_u - 2}$$

$$= (1+t)^{2r-c_1-1+i} \cdot \prod_{u=k+1}^{r} [(1+t) - m_u t]^{m_u - 2}$$

$$= (1+t)^{2r-c_1-1+i} \cdot \sum_{j=0}^{c_1-2r+k} \sum_{j_{k+1}+\ldots+j_r=j}$$

$$\cdot \prod_{u=k+1}^{r} \binom{m_u - 2}{j_u}(-m_u t)^{j_u} \cdot (1+t)^{m_u - 2 - j_u}$$

$$= \sum_{j=0}^{c_1-2r+k} \sum_{j_{k+1}+\ldots+j_r=j} \prod_{u=k+1}^{r} \binom{m_u - 2}{j_u}(-m_u t)^{j_u} \cdot (1+t)^{i+k-1-j}$$

since $\sum_{u=k+1}^{r} (m_u - 2 - j_u) = c_1 - 2r + k - j$. So $W_i$ is the coefficient of $t^i$ in

$$\prod_{u=k+1}^{r}(1 - m_u t)^{m_u - 2} \cdot \sum_{j=0}^{+\infty} \binom{j+k-1}{k-1} t^j = \prod_{u=k+1}^{r}(1 - m_u t)^{m_u - 2} \cdot \frac{1}{(1-t)^k}$$

$$= \prod_{u=1}^{r}(1 - m_u t)^{m_u - 2}. \quad \square$$

**Proposition 5.14.** *Let $V = \oplus_{i=1}^{r} \mathcal{O}_{\mathbb{P}^n}(m_i)$ where $m_i \geq 1$ for each $i$ and*

$$\sum_{i=1}^{r} m_i < \min(2r, (n + 1 + 2r)/2).$$

*Then the first relation $f_\omega^1$ for the quantum cohomology ring $H_\omega^*(\mathbb{P}(V); \mathbb{Z})$ is*

$$h^{n+1} = \prod_{u=1}^{r} (\xi - m_u h)^{m_u - 1} \cdot e^{-t(n+1+r-\sum_{i=1}^{r} m_i)}. \tag{5.15}$$

*Proof.* We may assume that $1 = m_1 = \ldots = m_k < m_{k+1} \leq \ldots \leq m_r$. Since the conclusion clearly holds when $k = r$, we also assume that $k < r$. Let $c_1 = \sum_{i=1}^{r} m_i$. Notice that the conditions in Corollary 4.7 (i) are satisfied. Thus,

$$h^{n+1} = \left( \sum_{i=0}^{c_1 - r} a_i \cdot h^i \cdot \xi^{c_1 - r - i} \right) \cdot e^{-t(n+1+r-c_1)}.$$

More directly, putting $\tilde{n} = \left[ \frac{n+1}{2} \right]$, then $\tilde{n} < -K_{\mathbb{P}(V)}(A_2) = (n+1+r-c_1)$, and $(n+1-\tilde{n}) < -K_{\mathbb{P}(V)}(A_2)$ unless $n$ is even and $c_1 = (n+2r)/2$. From the proofs in Theorem 4.1 and Corollary 4.7 (i) for the first relation $f_\omega^1$, we have $h^{\tilde{n}} = h_{\tilde{n}}$, and $h^{n+1-\tilde{n}} = h_{n+1-\tilde{n}}$ unless $n$ is even and $c_1 = (n+2r)/2$. Moreover, if $n$ is even and $c_1 = (n+2r)/2$, then $h^{n+1-\tilde{n}} = h \cdot h^{n-\tilde{n}} = h \cdot h_{n-\tilde{n}} = h_{n+1-\tilde{n}} + (h \cdot h_{n-\tilde{n}})_{A_2} \cdot e^{-t(n+1+r-c_1)}$. Since $(h \cdot h_{n-\tilde{n}})_{A_2}$ is of degree zero, $(h \cdot h_{n-\tilde{n}})_{A_2} = \Phi_{(A_2,0)}(h, h_{n-\tilde{n}}, h_n \xi_{r-1})$. Since $1 \leq k < r$, we can choose a point $q_0$ in $\mathbb{P}(V)$ representing the homology class $(h_n \xi_{r-1})_*$ such that the point $q_0$ is not contained in the $(k-1)$-dimensional linear subspace

$$\mathbb{P}^{k-1} = \mathbb{P}((\mathcal{O}_{\mathbb{P}^n}(1)^{\oplus k})|_{\pi(q_0)}) \subset \mathbb{P}(V|_{\pi(q_0)}) \cong \mathbb{P}^{r-1}.$$

Note that for every $f \in \mathfrak{M}(A_2, 0)$, $\text{Im}(f) = \ell \times \{q\}$ for some line $\ell \subset \mathbb{P}^n$ and some point $q \in \mathbb{P}^{k-1}$. Thus $\text{Im}(f)$ can not pass $q_0$. As in the proof of Lemma 3.7, we conclude that $\Phi_{(A_2,0)}(h, h_{n-\tilde{n}}, h_n \xi_{r-1}) = 0$. Therefore, $h^{n+1-\tilde{n}} = h_{n+1-\tilde{n}}$. So

$$h^{n+1} = h^{\tilde{n}} \cdot h^{n+1-\tilde{n}} = h_{\tilde{n}} \cdot h_{n+1-\tilde{n}}.$$

By similar arguments in the proofs of Theorem 4.1 and Corollary 4.7 (i) for the first relation $f_\omega^1$, we see that if $(h_{\tilde{n}} \cdot h_{n+1-\tilde{n}})_A \neq 0$, then $A = 0, A_2$. Thus

$$h^{n+1} = h_{n+1} + (h_{\tilde{n}} \cdot h_{n+1-\tilde{n}})_{A_2} \cdot e^{-t(n+1+r-c_1)} = (h_{\tilde{n}} \cdot h_{n+1-\tilde{n}})_{A_2} \cdot e^{-t(n+1+r-c_1)}.$$

So it suffices to show that $(h_{\tilde{n}} \cdot h_{n+1-\tilde{n}})_{A_2} = \prod_{u=1}^{r} (\xi - m_u h)^{m_u - 1}$. Note that

$$\prod_{u=1}^{r} (\xi - m_u h)^{m_u - 1} = \prod_{u=1}^{r} (\xi - m_u h)_{m_u - 1}$$

where the right-hand-side stands for the product in the ordinary cohomology. Thus we need to show that $(h_{\tilde{n}} \cdot h_{n+1-\tilde{n}})_{A_2} = \prod_{u=1}^{r} (\xi - m_u h)_{m_u - 1}$, or equivalently,

$$\Phi_{(A_2,0)}(h_{\tilde{n}}, h_{n+1-\tilde{n}}, h_{n-i} \xi_{2r-c_1-1+i}) = \prod_{u=1}^{r} (\xi - m_u h)_{m_u - 1} h_{n-i} \xi_{2r-c_1-1+i} \quad (5.16)$$

for $0 \leq i \leq (c_1 - i)$. The left-hand-side of (5.16) is computed in Lemma 5.12.

Denote the right-hand-side of (5.16) by $\tilde{W}_i$. Let $s_i$ be the $i$-th Segre class of $V$. Then we have $s_i = (-1)^i \cdot \sum_{j_1 + \ldots + j_r = i} \prod_{u=1}^{r} m_u^{j_u}$ and

$$\sum_{i=0}^{+\infty} (-1)^i s_i t^i = \prod_{u=1}^{r} \frac{1}{1 - m_u t}. \quad (5.17)$$

Moreover from the second relation in (1.1), we obtain for $i \geq r$,

$$\xi_i = (-1)^{i-(r-1)} s_{i-(r-1)} \xi_{r-1} + \text{ (terms with exponentials of } \xi \text{ less than } (r-1)).$$

It follows from the right-hand-side of (5.16) that $\tilde{W}_i$ is equal to

$$\sum_{j=0}^{c_1-r} \sum_{j_1+\ldots+j_r=j} \prod_{u=1}^{r} \binom{m_u-1}{j_u} \xi_{m_u-1-j_u} (-m_u h)_{j_u} h_{n-i} \xi_{2r-c_1-1+i}$$

$$= \sum_{j=0}^{i} \sum_{j_1+\ldots+j_r=j} \prod_{u=1}^{r} \binom{m_u-1}{j_u} (-m_u)^{j_u} h_{n-i+j} \xi_{r-1+i-j}$$

$$= \sum_{j=0}^{i} (-1)^{i-j} s_{i-j} \sum_{j_1+\ldots+j_r=j} \prod_{u=1}^{r} \binom{m_u-1}{j_u} (-m_u)^{j_u}.$$

Therefore, the formal power series $\sum_{i=0}^{+\infty} \tilde{W}_i t^i$ is equal to

$$\sum_{i=0}^{+\infty} \sum_{j=0}^{i} (-1)^{i-j} s_{i-j} t^{i-j} \sum_{j_1+\ldots+j_r=j} \prod_{u=1}^{r} \binom{m_u-1}{j_u} (-m_u t)^{j_u}$$

$$= \sum_{j=0}^{+\infty} \sum_{i=j}^{+\infty} (-1)^{i-j} s_{i-j} t^{i-j} \sum_{j_1+\ldots+j_r=j} \prod_{u=1}^{r} \binom{m_u-1}{j_u} (-m_u t)^{j_u}$$

$$= \sum_{j=0}^{+\infty} \sum_{i=0}^{+\infty} (-1)^{i} s_i t^i \sum_{j_1+\ldots+j_r=j} \prod_{u=1}^{r} \binom{m_u-1}{j_u} (-m_u t)^{j_u}$$

$$= \sum_{j=0}^{+\infty} \prod_{u=1}^{r} \frac{1}{1-m_u t} \sum_{j_1+\ldots+j_r=j} \prod_{u=1}^{r} \binom{m_u-1}{j_u} (-m_u t)^{j_u}$$

$$= \prod_{u=1}^{r} \frac{1}{1-m_u t} \sum_{j=0}^{+\infty} \sum_{j_1+\ldots+j_r=j} \prod_{u=1}^{r} \binom{m_u-1}{j_u} (-m_u t)^{j_u}$$

$$= \prod_{u=1}^{r} \frac{1}{1-m_u t} \prod_{u=1}^{r} (1-m_u t)^{m_u-1}$$

$$= \prod_{u=1}^{r} (1-m_u t)^{m_u-2}$$

where we have applied (5.17) in the third equality. By Lemma 5.12, $\tilde{W}_i = W_i$ for $0 \leq i \leq (c_1 - r)$. Hence the formule (5.16) and (5.15) hold. $\square$

It turns out that under certain conditions on the integers $m_i$, the second relation $f_\omega^2$ for the quantum cohomology ring $H_\omega^*(\mathbb{P}(V); \mathbb{Z})$ is much easier to be determined. Note that the second relation $f^2$ in (1.1) can be rewritten as

$$\prod_{i=1}^{r} (\xi - m_i h) = 0 \tag{5.18}$$

where the left-hand-side stands for the product in the ordinary cohomology ring.

**Proposition 5.19.** *Let $V = \oplus_{i=1}^{r} \mathcal{O}_{\mathbb{P}^n}(m_i)$ where $m_i \geq 1$ for each $i$, $m_i = 1$ for some $i$, and $\sum_{i=1}^{r} m_i < (2n+2+r)/2$. Then the second relation $f_\omega^2$ for the quantum cohomology ring $H_\omega^*(\mathbb{P}(V); \mathbb{Z})$ is*

$$\prod_{i=1}^{r}(\xi - m_i h) = e^{-tr} \tag{5.20}$$

*where the left-hand-side stands for the product in the quantum cohomology ring.*

*Proof.* We may assume that $1 = m_1 = \ldots = m_k < m_{k+1} \leq \ldots \leq m_r$. So $k \geq 1$. We notice that the conditions in Corollary 4.7 (ii) are satisfied. From the proofs of Theorem 4.1 (i) and Corollary 4.7 (ii), we see that the quantum corrections to the second relation (5.18) can only come from the classes $A_1, A_2$; moreover, the quantum correction from $A_1$ is $e^{-tr}$. Thus it suffices to show that the quantum correction from $A_2$ is zero. In view of (3.3), it suffices to show that

$$\Phi_{(A_2,0)}(\xi - m_1 h, \ldots, \xi - m_r h, \alpha) = 0$$

for every $\alpha \in H^*(\mathbb{P}(V); \mathbb{Z})$. For $1 \leq i \leq r$, let $V_i$ be the subbundle of $V$:

$$V_i = \mathcal{O}_{\mathbb{P}^n}(m_1) \oplus \ldots \oplus \mathcal{O}_{\mathbb{P}^n}(m_{i-1}) \oplus \mathcal{O}_{\mathbb{P}^n}(m_{i+1}) \oplus \ldots \oplus \mathcal{O}_{\mathbb{P}^n}(m_r),$$

and let $B_i = \mathbb{P}(V_i)$ be the codimension-1 subvariety of $\mathbb{P}(V)$ induced by the projection $V \to V_i \to 0$. Then the fundamental class of $B_i$ is $(\xi - m_i h)$. As in the proof of Lemma 3.7, we need only to show that if $f \in \mathfrak{M}(A_2, 0)$, then the image $\text{Im}(f)$ can not intersect with $B_1, \ldots, B_r$ simultaneously. In fact, we will show that $\text{Im}(f)$ can not intersect with $B_1, \ldots, B_k$ simultaneously. Indeed, $\text{Im}(f)$ is of the form

$$\text{Im}(f) = \ell \times \{q\} \subset \ell \times \mathbb{P}^{k-1} \subset \mathbb{P}^n \times \mathbb{P}^{k-1} = \mathbb{P}(\mathcal{O}_{\mathbb{P}^n}(1)^{\oplus k})$$

for some line $\ell \subset \mathbb{P}^n$, and $B_i|_{\pi^{-1}(\ell)} = \mathbb{P}(V_i|_\ell)$. Put $p = \pi(q) \in \mathbb{P}^n$, and

$$V|_p = \oplus_{i=1}^{k} \mathbb{C} \cdot e_i \oplus (\oplus_{i=k+1}^{r} \mathcal{O}_{\mathbb{P}^n}(m_i)|_p)$$

where $e_i$ is a global section of $\mathcal{O}_{\mathbb{P}^n}(m_i) = \mathcal{O}_{\mathbb{P}^n}(1)$ for $i \leq k$. Now the point $q$ is identified with $\mathbb{C} \cdot v$ for some nonzero vector $v \in \oplus_{i=1}^{k} \mathbb{C} \cdot e_i$. Let $v = \sum_{i=1}^{k} a_i e_i$. Since $\ell \times \{q\}$ and $B_i$ ($1 \leq i \leq k$) intersect, the one-dimensional vector space $\mathbb{C} \cdot v$ is also contained in $(V_i)|_p$. It follows that $a_i = 0$ for every $i$ with $1 \leq i \leq k$. But this is impossible since $v$ is a nonzero vector. $\square$

In summary, we partially verify Batyrev's conjecture.

**Theorem 5.21.** *Let $V = \oplus_{i=1}^{r} \mathcal{O}_{\mathbb{P}^n}(m_i)$ where $m_i \geq 1$ for each $i$ and*

$$\sum_{i=1}^{r} m_i < \min(2r, (n+1+2r)/2, (2n+2+r)/2).$$

*Then the quantum cohomology $H_\omega^*(\mathbb{P}(V); \mathbb{Z})$ is generated by $h$ and $\xi$ with relations*

$$h^{n+1} = \prod_{i=1}^{r}(\xi - m_i h)^{m_i - 1} \cdot e^{-t(n+1+r-\sum_{i=1}^{r} m_i)} \qquad \text{and} \qquad \prod_{i=1}^{r}(\xi - m_i h) = e^{-tr}.$$

*Proof.* Follows immediately from Propositions 5.14 and 5.19. $\square$

## 6. Examples

In this section, we shall determine the quantum cohomology of $\mathbb{P}(V)$ for ample bundles $V$ over $\mathbb{P}^n$ with $2 \le r \le n$ and $c_1 = r+1$. In these cases, $V|_\ell = \mathcal{O}_\ell(1)^{\oplus(r-1)} \oplus \mathcal{O}_\ell(2)$ for every line $\ell \subset \mathbb{P}^n$. In particular, $V$ is a uniform bundle. If $r < n$, then by the Theorem 3.2.3 in [10], $V = \mathcal{O}_{\mathbb{P}^n}(1)^{\oplus(r-1)} \oplus \mathcal{O}_{\mathbb{P}^n}(2)$; if $r = n$, then by the results on pp.71-72 in [10], $V = \mathcal{O}_{\mathbb{P}^n}(1)^{\oplus(n-1)} \oplus \mathcal{O}_{\mathbb{P}^n}(2)$ or $V = T_{\mathbb{P}^n}$ the tangent bundle of $\mathbb{P}^n$. When $V = \mathcal{O}_{\mathbb{P}^n}(1)^{\oplus(r-1)} \oplus \mathcal{O}_{\mathbb{P}^n}(2)$ with $r \le n$, the conditions in Theorem 5.21 are satisfied, so the quantum cohomology ring $H^*_\omega(\mathbb{P}(V);\mathbb{Z})$ is the ring generated by $h$ and $\xi$ with two relations

$$h^{n+1} = (\xi - 2h) \cdot e^{-t(n+1+r-c_1)} \quad \text{and} \quad (\xi - h)^{r-1}(\xi - 2h) = e^{-tr}.$$

In the rest of this section, we compute the quantum cohomology of $\mathbb{P}(T_{\mathbb{P}^n})$. It is well-known that $(\xi - h)$ is a nef divisor on $\mathbb{P}(T_{\mathbb{P}^n})$, and the two extremal rays of $\mathbb{P}(T_{\mathbb{P}^n})$ are generated by $A_1 = (h_n \xi_{n-2})_*$ and $A_2 = (h_{n-1}\xi_{n-1} - nh_n\xi_{n-2})_*$. Moreover, $A_2$ is represented by smooth rational curves in $\mathbb{P}(T_{\mathbb{P}^n})$ induced by the surjective maps $T_{\mathbb{P}^n}|_\ell \to \mathcal{O}_\ell(1) \to 0$ for lines $\ell \subset \mathbb{P}^n$. Since $c_1 = n+1$ and $n \ge 2$, the assumptions in Corollary 4.7 are satisfied, so the quantum cohomology ring $H^*_\omega(\mathbb{P}(T_{\mathbb{P}^n});\mathbb{Z})$ is the ring generated by $h$ and $\xi$ with two relations

$$h^{n+1} = (a_1 h + \xi) \cdot e^{-tn} \quad \text{and} \quad \sum_{i=0}^{n}(-1)^i c_i \cdot h^i \cdot \xi^{n-i} = (1+b_0) \cdot e^{-tn}. \qquad (6.1)$$

More precisely, putting $H'_2 = H_2(\mathbb{P}(V);\mathbb{Z}) - \mathbb{Z} \cdot (h_n \xi_{n-2})_*$, then we see from the proof of Corollary 4.7 (i) that the only homology class $A \in H'_2$ which has nonzero contributions to the quantum corrections in (4.4) is $A = A_2$. Thus by (4.4),

$$h \cdot h_p = \begin{cases} h_{p+1}, & \text{if } p \le n-2 \\ h_n + a'_1 \cdot e^{-tn}, & \text{if } p = n-1 \\ h_{n+1} + (a'_2 h + a'_3 \xi) \cdot e^{-tn}, & \text{if } p = n. \end{cases} \qquad (6.2)$$

where $a'_1 = \Phi_{(A_2,0)}(h, h_{n-1}, h_n \xi_{n-1})$, $a'_3 = \Phi_{(A_2,0)}(h, h_n, h_n \xi_{n-2})$, and

$$a'_2 = \Phi_{(A_2,0)}(h, h_n, h_{n-1}\xi_{n-1}) - c_1 a'_3.$$

By Lemma 3.8, $a'_3 = 1$. Thus $a_1 = (a'_1 + a'_2)$ and the first relation $f^1_\omega$ in (6.1) is

$$h^{n+1} = ((a'_1 + a'_2)h + \xi) \cdot e^{-tn} \qquad (6.3)$$

Similarly, from the proof of Corollary 4.7 (ii), we see that the only homology class $A \in H'_2$ which has nonzero contributions to the quantum corrections in (4.5) and (4.6) is also $A = A_2$. By (4.5), $\xi \cdot \xi_p = \xi_{p+1}$ if $p < n-1$, and $\xi \cdot \xi_{n-1} = \xi_n + e^{-tn} + b_2^{(n)} \cdot e^{-tn}$ where $b_2^{(n)} = \Phi_{(A_2,0)}(\xi, \xi_{n-1}, h_n \xi_{n-1})$. Thus,

$$\xi^p = \begin{cases} \xi_p, & \text{if } p < n \\ \xi_n + (1 + b_2^{(n)}) \cdot e^{-tn}, & \text{if } p = n \end{cases} \qquad (6.4)$$

By (6.2), we have $h \cdot h_p = h_{p+1}$ if $p < n-1$, and $h \cdot h_{n-1} = h_n + b_2^{(0)} \cdot e^{-tn}$ where $b_2^{(0)} = a_1' = \Phi_{(A_2,0)}(h_{n-1}, h, h_n\xi_{n-1})$. Thus, we obtain

$$h^p = \begin{cases} h_p, & \text{if } p < n \\ h_n + b_2^{(0)} \cdot e^{-tn}, & \text{if } p = n \end{cases} \quad (6.5)$$

By (4.6), for $1 \leq i \leq (n-1)$, $h_{n-i} \cdot \xi_i = h_{n-i}\xi_i + b_2^{(i)} \cdot e^{-tn}$ where $b_2^{(i)} = \Phi_{(A_2,0)}(h_{n-i}, \xi_i, h_n\xi_{n-1})$. Thus by (6.4) and (6.5), we have

$$h^{n-i} \cdot \xi^i = h_{n-i} \cdot \xi_i = h_{n-i}\xi_i + b_2^{(i)} \cdot e^{-tn}. \quad (6.6)$$

Since $\sum_{i=0}^{n}(-1)^i c_i \cdot h_i\xi_{n-i} = 0$, it follows from (6.4), (6.5), (6.6) that

$$\sum_{i=0}^{n}(-1)^i c_i \cdot h^i \cdot \xi^{n-i} = (1 + \sum_{i=0}^{n}(-1)^i c_i b_2^{(n-i)}) \cdot e^{-tn}. \quad (6.7)$$

Next, we compute the above integers $a_1', a_2'$, and $b_2^{(i)}$ where $0 \leq i \leq n$.

**Lemma 6.8.** Let $V = T_{\mathbb{P}^n}$ with $n \geq 2$ and $A_2 = (h_{n-1}\xi_{n-1} - nh_n\xi_{n-2})_*$.
  (i) $\Phi_{(A_2,0)}(h, h_n, h_{n-1}\xi_{n-1}) = n$;
  (ii) Let $\alpha = h_j\xi_k$ and $\beta = h_s\xi_t$ where $j, k, s, t$ are nonnegative integers such that $\max(j,k) > 0$, $\max(s,t) > 0$, and $(j+k+s+t) = n$. Then,

$$\Phi_{(A_2,0)}(\alpha, \beta, h_n\xi_{n-1}) = 1.$$

*Proof.* (i) By Lemma 2.2 (iv), $\mathfrak{M}(A_2,0)/PSL(2;\mathbb{C})$ is compact. By (3.17), we have $h^1(N_f) = 0$ for every $f \in \mathfrak{M}(A_2,0)$. Thus, $\mathfrak{M}(A_2,0)/PSL(2;\mathbb{C})$ is also smooth. Fix a line $\ell_0$ in $\mathbb{P}^n$. Let $g: \ell_0 \to \mathbb{P}(T_{\mathbb{P}^n}|_{\ell_0}) \subset \mathbb{P}(T_{\mathbb{P}^n})$ be the embedding induced by the natural projection $T_{\mathbb{P}^n}|_{\ell_0} = \mathcal{O}_{\ell_0}(1)^{\oplus(n-1)} \oplus \mathcal{O}_{\ell_0}(2) \to \mathcal{O}_{\ell_0}(2) \to 0$. Since $h([g(\ell_0)]) = 1$ and $\xi([g(\ell_0)]) = 2$, we have $[g(\ell_0)] = [h_{n-1}\xi_{n-1} - (n-1)h_n\xi_{n-2}]_*$. So $h_{n-1}\xi_{n-1} = [g(\ell_0)]_* + (n-1)h_n\xi_{n-2}$, and

$$\Phi_{(A_2,0)}(h, h_n, h_{n-1}\xi_{n-1}) = \Phi_{(A_2,0)}(h, h_n, [g(\ell_0)]_*) + (n-1)\Phi_{(A_2,0)}(h, h_n, h_n\xi_{n-2}).$$

By Lemma 3.8, it suffices to show that $\Phi_{(A_2,0)}(h, h_n, [g(\ell_0)]_*) = 1$. Let $B$ and $C$ be the subvarieties of $\mathbb{P}(T_{\mathbb{P}^n})$ in general position, whose homology classes are Poincaré dual to $h$ and $h_n$ respectively. Then the homology classes of $\pi(B)$ and $\pi(C)$ in $\mathbb{P}^n$ are Poincaré dual to $h$ and $h_n$ respectively. Let $f: \ell \to \mathbb{P}(T_{\mathbb{P}^n})$ be a morphism in $\mathfrak{M}(A_2,0)$ induced by a surjective map $T_{\mathbb{P}^n}|_\ell \to \mathcal{O}_\ell(1) \to 0$ for some line $\ell \subset \mathbb{P}^n$. If the image $\text{Im}(f)$ intersects with $B, C$, and $g(\ell_0)$, then $\ell$ intersects with $\pi(B)$, $\pi(C)$, and $\pi(g(\ell_0)) = \ell_0$. In other words, $\ell$ passes through the point $\pi(C)$ and intersects with $\ell_0$. Moreover, putting $p = \ell \cap \ell_0$ and noticing that every surjective map $T_{\mathbb{P}^n}|_\ell \to \mathcal{O}_\ell(1) \to 0$ factors through the natural projection $T_{\mathbb{P}^n}|_\ell = \mathcal{O}_\ell(1)^{(n-1)} \oplus \mathcal{O}_\ell(2) \to \mathcal{O}_\ell(1)^{(n-1)}$, we conclude that the $(n-1)$-dimensional

subspace $(\mathcal{O}_\ell(1)^{(n-1)})|_p$ in $(T_{\mathbb{P}^n}|_\ell)|_p = T_{p,\mathbb{P}^n}$ must contain the 1-dimensional subspace $(\mathcal{O}_{\ell_0}(2))|_p$ in $(T_{\mathbb{P}^n}|_{\ell_0})|_p = T_{p,\mathbb{P}^n}$. Conversely, let $p \in \ell_0$ and let $\ell_p$ be the unique line connecting the two points $\pi(C)$ and $p$. If the $(n-1)$-dimensional subspace $(\mathcal{O}_{\ell_p}(1)^{(n-1)})|_p$ in $(T_{\mathbb{P}^n}|_{\ell_p})|_p = T_{p,\mathbb{P}^n}$ contains the 1-dimensional subspace $(\mathcal{O}_{\ell_0}(2))|_p$ in $(T_{\mathbb{P}^n}|_{\ell_0})|_p = T_{p,\mathbb{P}^n}$, then there exists a unique surjective map $T_{\mathbb{P}^n}|_{\ell_p} \to \mathcal{O}_{\ell_p}(1) \to 0$ such that the image of the induced morphism $f: \ell_p \to \mathbb{P}(T_{\mathbb{P}^n})$ intersects $g(\ell_0)$ at the point $g(p)$. Since there exists a unique point $p \in \ell_0$ such that the $(n-1)$-dimensional subspace $(\mathcal{O}_{\ell_p}(1)^{(n-1)})|_p$ in $(T_{\mathbb{P}^n}|_{\ell_p})|_p = T_{p,\mathbb{P}^n}$ contains the 1-dimensional subspace $(\mathcal{O}_{\ell_0}(2))|_p$ in $(T_{\mathbb{P}^n}|_{\ell_0})|_p = T_{p,\mathbb{P}^n}$, it follows that

$$\Phi_{(A_2,0)}(h, h_n, [g(\ell_0)]_*) = 1.$$

(ii) It is well-known (see p.176 of [7]) that there is an exact sequence

$$0 \to \mathcal{O}_{\mathbb{P}^n} \to \mathcal{O}_{\mathbb{P}^n}(1)^{\oplus(n+1)} \to T_{\mathbb{P}^n} \to 0. \qquad (6.9)$$

The surjective map $\mathcal{O}_{\mathbb{P}^n}(1)^{\oplus(n+1)} \to T_{\mathbb{P}^n} \to 0$ induces the inclusion $\phi: \mathbb{P}(T_{\mathbb{P}^n}) \subset \mathbb{P}^n \times \mathbb{P}^n$ such that $\xi$ is the restriction of the $(1,1)$ class in $\mathbb{P}^n \times \mathbb{P}^n$. Let $B, C, q_0$ be the subvarieties of $\mathbb{P}(T_{\mathbb{P}^n})$ in general position, whose homology classes are Poincaré dual to $\alpha, \beta, h_n \xi_{n-1}$ respectively. Then $q_0$ is a point. Put $p_0 = \pi(q_0) \in \mathbb{P}^n$. Now the morphisms in $\mathfrak{M}(A_2, 0)$ are of the forms $f: \ell \to \mathbb{P}(T_{\mathbb{P}^n})$ induced by surjective maps $T_{\mathbb{P}^n}|_\ell \to \mathcal{O}_\ell(1) \to 0$ for lines $\ell \subset \mathbb{P}^n$. If the image $\mathrm{Im}(f)$ passes $q_0$, then the line $\ell$ passes $p_0$ and $q_0$ is contained in the hyperplane

$$\mathbb{P}^{n-2} = \mathbb{P}((\mathcal{O}_\ell(1)^{\oplus(n-1)})|_{p_0}) \subset \mathbb{P}((T_{\mathbb{P}^n}|_\ell)|_{p_0}) = \pi^{-1}(p_0) = \mathbb{P}^{n-1}.$$

Conversely, if $\ell$ passes $p_0$ and $q_0$ is contained in the hyperplane

$$\mathbb{P}^{n-2} = \mathbb{P}((\mathcal{O}_\ell(1)^{\oplus(n-1)})|_{p_0}) \subset \mathbb{P}((T_{\mathbb{P}^n}|_\ell)|_{p_0}) = \pi^{-1}(p_0) = \mathbb{P}^{n-1}, \qquad (6.10)$$

then there exists a unique $f \in \mathfrak{M}(A_2, 0)$ of the form $f: \ell \to \mathbb{P}(T_{\mathbb{P}^n})$ such that $\mathrm{Im}(f)$ passes $q_0$; moreover, putting $q_0 = (p_0, p_0') \in \mathbb{P}^n \times \mathbb{P}^n$ such that $\pi$ is the first projection of $\mathbb{P}^n \times \mathbb{P}^n$, then $\mathrm{Im}(f) = \ell \times \{p_0'\} \subset \mathbb{P}^n \times \mathbb{P}^n$. The set of all lines $\ell$ passing $p_0$ such that $q_0$ is contained in the hyperplane (6.10) is parameterized by an $(n-2)$-dimensional linear subspace $\mathbb{P}^{n-2}$ in $\mathbb{P}^n$ (the first factor in $\mathbb{P}^n \times \mathbb{P}^n$). It follows that the images $\mathrm{Im}(f) \subset \mathbb{P}(T_{\mathbb{P}^n})$ sweep a hyperplane

$$H \stackrel{\mathrm{def}}{=} \mathbb{P}^{n-1} \times \{p_0'\} \subset \mathbb{P}^n \times \{p_0'\}. \qquad (6.11)$$

Since $\xi$ is the restriction of the $(1,1)$ class in $\mathbb{P}^n \times \mathbb{P}^n$, $\xi|_H$ is the hyperplane class $\tilde{h}$ in $H = \mathbb{P}^{n-1} \times \{p_0'\} \cong \mathbb{P}^{n-1}$. Thus $\alpha|_H = \tilde{h}_{j+k}$ and $\beta|_H = \tilde{h}_{s+t}$. Since $(j + k + s + t) = n$ and $B$ and $C$ are in general position, there is a unique line in $H$ passing $q_0 = (p_0, p_0')$ and intersecting with $B$ and $C$. Therefore,

$$\Phi_{(A_2,0)}(\alpha, \beta, h_n \xi_{n-1}) = 1. \quad \square$$

Finally, we summarize the above computations and prove the following.

**Proposition 6.12.** *The quantum cohomology ring $H^*_\omega(\mathbb{P}(T_{\mathbb{P}^n}); \mathbb{Z})$ with $n \geq 2$ is the ring generated by $h$ and $\xi$ with the two relations:*

$$h^{n+1} = \xi \cdot e^{-tn} \qquad \text{and} \qquad \sum_{i=0}^{n}(-1)^i c_i \cdot h^i \cdot \xi^{n-i} = (1+(-1)^n) \cdot e^{-tn}.$$

*Proof.* By Lemma 6.8 (ii), $a'_1 = 1$. By Lemma 3.8, $a'_3 = 1$. By Lemma 6.8 (i),

$$a'_2 = \Phi_{(A_2,0)}(h, h_n, h_{n-1}\xi_{n-1}) - c_1 a'_3 = -1.$$

Thus by (6.3), the first relation $f^1_\omega$ is $h^{n+1} = \xi \cdot e^{-tn}$. By Lemma 6.8 (ii), $b_2^{(i)} = 1$ for $0 \leq i \leq n$. By (6.7), the second relation $f^2_\omega$ is $\sum_{i=0}^{n}(-1)^i c_i \cdot h^i \cdot \xi^{n-i} = (1 + \sum_{i=0}^{n}(-1)^i c_i) \cdot e^{-tn}$. From the exact sequence (6.9), $c_i = \binom{n+1}{i}$ for $0 \leq i \leq n$. Therefore, the relation $f^2_\omega$ is $\sum_{i=0}^{n}(-1)^i c_i \cdot h^i \cdot \xi^{n-i} = (1+(-1)^n) \cdot e^{-tn}$. □

Department of Mathematics, Oklahoma State University, Stillwater, OK 74078
*E-mail address*: zq@math.okstate.edu

Department of Mathematics, University of Utah, Salt Lake City, UT 84112
*E-mail address*: ruan@math.utah.edu